\documentclass{amsart}
\usepackage{hyperref}
\usepackage{graphicx,epsfig,subfigmat}
\usepackage{cite,amsmath,amssymb,epsfig,amsfonts,amscd,mathrsfs,enumerate,latexsym,url,tabularx,framed}
\usepackage{multirow,booktabs,algorithm,algorithmic,color,epstopdf}
\usepackage{chngcntr}
\counterwithin{section}{part}

\newcommand{\E}{\mathscr{E}}

\newcommand{\p}{\mathscr{P}}
\newcommand{\s}{\mathscr{S}}
\newcommand{\eg}{{\rm e.g., }}
\newcommand{\ie}{{\rm i.e., }}

\def\be {\begin{equation}}
\def\ee {\end{equation} }
\def\ba {\begin{eqnarray} }
\def\ea {\end{eqnarray} }
\def\bes {\begin{equation*} }
\def\ees {\end{equation*} }
\def\bas {\begin{eqnarray*} }
\def\eas {\end{eqnarray*} }
\def\bpr {\begin{proof} }
\def\epr {\end{proof} }

\newtheorem{theorem}{Theorem}[section]

\theoremstyle{definition}

\newtheorem{example}[theorem]{Example}

\newtheorem{remark}[theorem]{Remark}

\theoremstyle{question}

\newtheorem{defn}[theorem]{Definition}
\numberwithin{equation}{section}

\begin{document}

%\noindent {\footnotesize -- August 31, 2015.}\\[.5in]

\title[A user guide for {\tt Singularity}]{A user guide for {\tt Singularity}}

\author{Majid Gazor and Mahsa Kazemi}

\address{Department of Mathematical Sciences, Isfahan University of Technology, Isfahan 84156--83111, Iran.}

%\address{School of Mathematics, Institute for Research in Fundamental Sciences (IPM), P.O. Box: 19395--5746, Tehran, Iran.}

\email{mgazor@cc.iut.ac.ir; mahsa.kazemi@math.iut.ac.ir}



%\thanks{}

\begin{abstract}
This is a user guide for the first version of our developed Maple library, named {\tt Singularity}. The first version here is designed for the qualitative study of local real zeros of scalar smooth maps. This library will be extended for symbolic bifurcation analysis and control of different singularities including autonomous differential singular systems and local real zeros of multidimensional smooth maps. Many tools and techniques from computational algebraic geometry have been used to develop {\tt Singularity}. However, we here skip any reference on how this library is developed. This package is useful for both pedagogical and research purposes. {\tt Singularity} will be updated as our research progresses.
% and will be released for public access once our draft paper \cite{GazorKazemi} is peer-reviewed in a refereed journal.
\end{abstract}

\maketitle

\tableofcontents

\part{How to use {\tt Singularity}}
%There are two options on how to use {\tt Singularity}  with {\sc Maple} on your computer.
%The first is to use a ``read'' command and the second option is to install {\tt Singularity} as a sub-package of {\sc Maple} on your computer.
In order to use {\tt Singularity} in Windows and Linux platforms follow the instructions below accordingly.

\subsection{Read option}

The following command reads {\tt Singularity} through which all of the commands in this user-guide will be accessible.
\begin{itemize}
\item Windows users: First save \verb"SingularityLibrary.mpl" with an address like
\verb"C:\your address\SingularityLibrary.mpl". Next, open a {\it Maple} worksheet and type
\begin{verbatim}
read "C:\\your address\\SingularityLibrary.mpl".
\end{verbatim}
\item Linux users: Save \verb"SingularityLibrary.mpl" with an address like
{\tt /your address/SingularityLibrary.mpl}. Then, open a {\it Maple} worksheet and type
{\tt read "/your address/SingularityLibrary.mpl"}.
\end{itemize}

\subsection{Installing {\tt Singularity} on your personal computer}
You can install {\tt Singularity} as an additional sub-package of {\it Maple} on your personal computer through the following steps. Then, similar to built-in packages of {\it Maple}, you can use {\tt Singularity} after calling \verb"with(Singularity)".  The latter enlists all the commands in subsections \ref{ListCommands} and \ref{ListAlgCommands} as its output.

\begin{itemize}
\item Windows users: Type the following lines in a {\it Maple} worksheet

\verb"with(LibraryTools);"
\begin{verbatim}
Create("C:\\Program Files\\Maple 15\\lib\\Singularity.lib");
\end{verbatim}
\begin{verbatim}
read "C:\\your address\\SingularityModule.mpl".
\end{verbatim}
%be careful the module structure is at the top of SingularityLibrary.mpl which is called SingularityModule.mpl
\begin{verbatim}
savelib('Singularity',"C:\\Program Files\\Maple 15\\lib\\Singularity.lib");
\end{verbatim}
Then close the worksheet and open a new one.
\item Linux users: Save {\tt Singularity.mla} in a directory listed in the {\it value of the predefined variable libname}; see \url{https://www.maplesoft.com/support/help/maple/view.aspx?path=libname}.
\end{itemize}
%\verb"march"('create',"/home/singlib.mla");
%\verb"read""/home/modulesingularity.mpl";
%\verb"Save"(Singularity, singlib);
%\verb"with(Singularity)";
\part{Introduction}

In this user guide, we describe how to use the first version of our developed {\sc Maple} library (named {\tt Singularity}) for local bifurcation analysis of real zeros of scalar smooth maps; see \cite{GazorKazemi} for the main ideas. We remark that the term {\it singularity theory} has been used in many different Mathematics disciplines with essentially different objectives and tools but yet sometimes with similar terminologies; for examples of these see \cite{GazorKazemi}. For more detailed information, definitions, and related theorems in what we call here {\it singularity theory}, we refer the reader to \cite{Melbourne87,Keyfitz,GatermannLauterbach,GatermannHosten,Gaffney,MurdBook,GovaertsBook,GolubitskySchaefferBook}.

\section{Verification and warning note} \texttt{Singularity} is able to check and verify all of its computations. However, this sometimes adds an extra computational cost. This happens mainly for finding out the correct and suitable truncation degree and computational ring. Therefore, it is beneficial to skip the extra computations when it is not necessary. For an instance of benefit, consider that you need to obtain certain results for a large family of problems arising from the same origin. Therefore, you might be able to only check a few problems and conclude about the suitable truncation degree and computational ring for the whole family. Thereby, the commands of \texttt{Singularity} check and verify the output results unless it requires extra computation. In this case, a warning note of not verified output or possible errors is given; in these cases, a recommendation is always provided on how to verify or circumvent the problem. Lack of warning notes always indicates that the output results have been successfully verified.

%\subsection{Public access }

%We will release the first version of {\tt Singularity} for public access when our draft paper \cite{GazorKazemi} is
%accepted for publication in a refereed journal.

%\subsection{Built-in help }

%We have constructed a comprehensive built-in help for our library. It can be used similar to already existing help commands of {\sc Maple}.
%For example the command

%\verb"? Tangent Space"

%\noindent pops up a help-window explaining the tangent space for a singular germ and how it can be computed by using our library.

%\subsection{How to install?}

%In this section we will explain how one can install this package on Maple software running on his/her computer.
%This is to be completed, once our library is released for public access.

%However, you can alternatively save this library on a folder and use the following command in the command window of your Maple software:
%\verb"read" ``\verb"Folder's address\\Singularity.mw""
%\noindent Next, you will be able to use all available commands of our library.

\section{List of commands in singularity theory}\label{ListCommands}

The following list shows a complete list of commands from singularity theory that have so far been implemented in {\tt Singularity}:

\begin{itemize}
\item \verb"CheckSingularity"; section \ref{checksing},
\item \verb"Verify"; section \ref{SecVerify},
\item \verb"Normalform"; section \ref{SecNormalform},
\item \verb"UniversalUnfolding"; \verb"UniversalUnfoldingPars"; section \ref{SecUniversalUnfolding},
\item \verb"RecognitionProblem"; section \ref{SecRecognitionProblem},
\item \verb"CheckUniversal"; subsection \ref{SecCheckUniversal},
\item \verb"Transformation"; section \ref{SecTransformation},
\item \verb"TransitionSet"; subsection \ref{SecTransitionSet},
\item \verb"PersistentDiagram"; subsection \ref{SecPersistentDiagram},
\item \verb"NonPersistent"; subsection \ref{SecNonPersistent},
\item \verb"Intrinsic"; section \ref{SecIntrinsic},
\item \verb"AlgObjects", \verb"RT", \verb"T", \verb"P", \verb"S", \verb"TangentPerp", \verb"SPerp", \verb"IntrinsicGen"; section \ref{SecAlgObjects}. \\
\end{itemize}

\section{Commands from computational algebraic geometry} \label{ListAlgCommands}
The following enlists all implemented tools from computational algebraic geometry:

\begin{itemize}
  \item \verb"MultMatrix"; section \ref{SecMultMatrix},
  \item \verb"Division"; section \ref{SecDivision},
  \item \verb"StandardBasis"; section \ref{SecStandardBasis},
  \item \verb"ColonIdeal"; section \ref{SecColonIdeal},
  \item \verb"Normalset"; section \ref{SecNormalset}.\\
\end{itemize}

\part{Singularity theory }

The terminology ``singularity theory'' has been used to deal with many different problems in different mathematical disciplines. Singularity theory here refers the methodologies in dealing with the qualitative behavior of local zeros
\begin{eqnarray}
\label{eq:1}
g(x, \lambda )=0,
\end{eqnarray}
where \(x\in \mathbb{R}\) is a state variable and \(\lambda\) is a distinguished parameter. The cases of multi-dimensional parameter space are dealt with through the notions of unfolding. {\tt Singularity} will be soon enhanced to deal with the cases of multi-dimensional state variables.

In many real life problems at certain solutions, \(g(x, \lambda)\) is {\it singular}, \ie
 \begin{eqnarray*}
 g(x, \lambda)=g_{x}(x, \lambda)=0.
\end{eqnarray*}
A singular germ \(g(x, \lambda)\) subjected to smooth changes demonstrates surprising changes in the {\it qualitative properties} of the solutions, \eg changes in the number of solutions. This phenomenon is called a {\it bifurcation}.

\section{Qualitative properties }

We define the {\it qualitative properties} as the invariance of an equivalence relation.
The equivalence relation used in {\tt Singularity} is {\it contact equivalence} and is defined by
\begin{equation}
f\sim g \Leftrightarrow f(x, \lambda)= S(x, \lambda)g(X(x, \lambda), \Lambda(\lambda))
\end{equation}
where \(S(x, \lambda)>0\) while \((X, \Lambda)\) is locally a diffeomorphism such that \(X_{x}(x, \lambda)\) and \(\Lambda^{\prime}(\lambda)>0\).

\section{Check singularity}\label{checksing}
Consider the bifurcation problem \(G(x, \lambda, \alpha)=0, \alpha\in \mathbb{R}^m, m\in \mathbb{Z}^{\geq0}.\) 
The command \(\verb"CheckSingularity"\) computes the singular points and their associated parameter varieties. In particular for a bifurcation problem, this command is designed to compute the varieties in the parameters space and values for the state variable \(x\) where the bifurcation problem has a given singularity. \\

\begin{tabularx}{\textwidth}{l|X}
\textbf{Command/the default options} & \textbf{Description} \\
\hline
  \verb"CheckSingularity"(\(G\), \verb"Pars", \verb"Vars")& computes the singular points of \(G\).
\end{tabularx}\\

Here \( \verb"Pars"\) represents the list of parameters \(\alpha\) (but not the distinguished parameter). These parameters are usually called unfoldings, bifurcation parameters, controller inputs, or the perturbation parameters. \(\verb"Vars"\) stands for the state variable \(x\) and the distinguished parameter. This provides a freedom for the user to work with his/her own choice of these variables. 

\subsection{Options}

\begin{itemize}
\item \verb"'type'='h"\((x,\lambda)\)'; where \(h\) is a type of singularity. This option derives the variety in the parameter and state variable space
\verb"Pars" and \verb"Vars" at which \(G\) has the singularity of \(h\)-type. We remark that the singularity type \(h\) is ALWAYS assumed to be in terms of the state variable \(x\) and distinguished parameter \(\lambda\) and their associated base point as the origin (zero). However, the state variable and distinguished parameter associated with germ \(G\) are determined in the \(\verb'Vars'\)-list. Further, the base point of the germ \(G\) can be different from the origin via the option \verb"SingularPoint" described below.
\item \verb"ParametricVariety"; calculates the parametric variety over which \(G\) attains a singularity of \(h\)-type.
\item \verb"'interval'"; computes the parametric variety over which \(G\) contains the singularity \(G\) provided that \verb"Pars" may only vary in the specified interval.
\item \verb"'ParsPoint'"; computes the parametric variety over which \(G\) acquires the singularity of \(h\)-type when some of the parameters in \verb"Pars"
receive the input values in \verb"ParsPoint" option.
\item \verb"'VarsPoint'"; derives the parametric variety over which \(G\) contains the singularity of \(h\)-type when some of the variables in \verb"Vars"
receive the input values in \verb"VarsPoint" option.
\item \verb"'plot'"; plots the parametric variety over which \(G\) contains the \(h\)-type singularity provided that \verb"Vars"/\verb"Pars" receive a specified input value via the option \verb"'VarsPoint'"/\verb"'ParsPoint'".
\item \verb"'SingularPoint'"='\([a,b]\)'; the computations are performed when the base point is given as \(a\) for the state variable and \(b\) for the
distinguished parameter.
\item \verb"NonZero"; provides the necessary non-zero conditions required for the singular germ \(G\) to have a \(h(x, \lambda)\)-type singularity. We
remark that some provided output conditions in using this option appear as equalities accompanied with some auxiliary parameters rather than nonzero conditions, \ie inequality. However, these types of output can often be readily translated into nonzero conditions.
\end{itemize}

\begin{example}
The default command \verb"CheckSingularity"\((x^3-\sin(\lambda-1)x, [x, \lambda])\) gives \((0, 1)\), \ie the germ \(x^3-\sin(\lambda-1)x\) is singular at \((x, \lambda)= (0, 1)\).

Now we provide several examples on how to use different options as follows:

\begin{enumerate}
\item Consider the command \verb"CheckSingularity"\((x^3-\lambda x+\alpha+\beta x^2, [\alpha, \beta], [x, \lambda]\), \verb"'type'"='\(x^2-\lambda\)', \verb"ParametricVariety").
This returns

``The bifurcation problem has singularity of \verb"'type'"='\(x^2-\lambda\)' for all values of parameters.''

\item The command
\verb"CheckSingularity"\((y^3-\eta y+\alpha+\beta y^2, [\alpha, \beta], [y, \eta]\), \verb"'type'"='\(x^3+\lambda^2\)', \verb"ParametricVariety") leads to

`` There is no singularity of this type in this bifurcation problem.''

\item The command
\verb"CheckSingularity"\((x^3+b\sin(\lambda^2+\alpha \exp(\lambda x)-\alpha)+\beta^2 x\lambda, [b, \alpha, \beta], [x, \lambda]\), \verb"'SingularPoint'"='\([0,0]\)', \verb"'type'"='\(x^3+\lambda^2\)', \verb"ParametricVariety") gives rise to

``The bifurcation problem has singularity \verb"'type'"='\(x^3+\lambda^2\)' when \(\lbrace \alpha b+\beta^2=0, b\neq 0\rbrace\).''
  \item The command \verb"CheckSingularity"\((z^3-\theta z+\alpha+\beta z^2, [\alpha, \beta], [z, \theta]\), \verb"'ParsPoint'"='\([\alpha=0]\)',\verb"'type'"='\(x^2-\lambda\)', \verb"ParametricVariety") concludes

``The bifurcation problem has singularity \verb"'type'"='\(x^2-\lambda\)' for all values of parameters.''
  \item \verb"CheckSingularity"\((x^3-\lambda x+\alpha+\beta x^2, [\alpha, \beta], [x, \lambda]\), \verb"'VarsPoint'"='\([\lambda=0]\)',\verb"'type'"='\(x^2-\lambda\)', \verb"ParametricVariety") returns

``The bifurcation problem has singularity \verb"'type'"='\(x^2-\lambda\)' when \(\lbrace 4\beta^3+27\alpha=0\rbrace\).''

\item \verb"CheckSingularity"\((z^3-\theta z+\alpha+\beta z^2, [\alpha, \beta], [z, \theta]\), \verb"'VarsPoint'"='\([\theta=0]\)',\verb"'type'"='\(x^2-\lambda\)', \verb"ParametricVariety", \verb"plot") gives rise to Figure \ref{plot_varspoint}.

\begin{figure}[h]
\begin{center}
\subfigure[\label{plot_varspoint}]{\includegraphics[width=.40\columnwidth,height=.4\columnwidth]{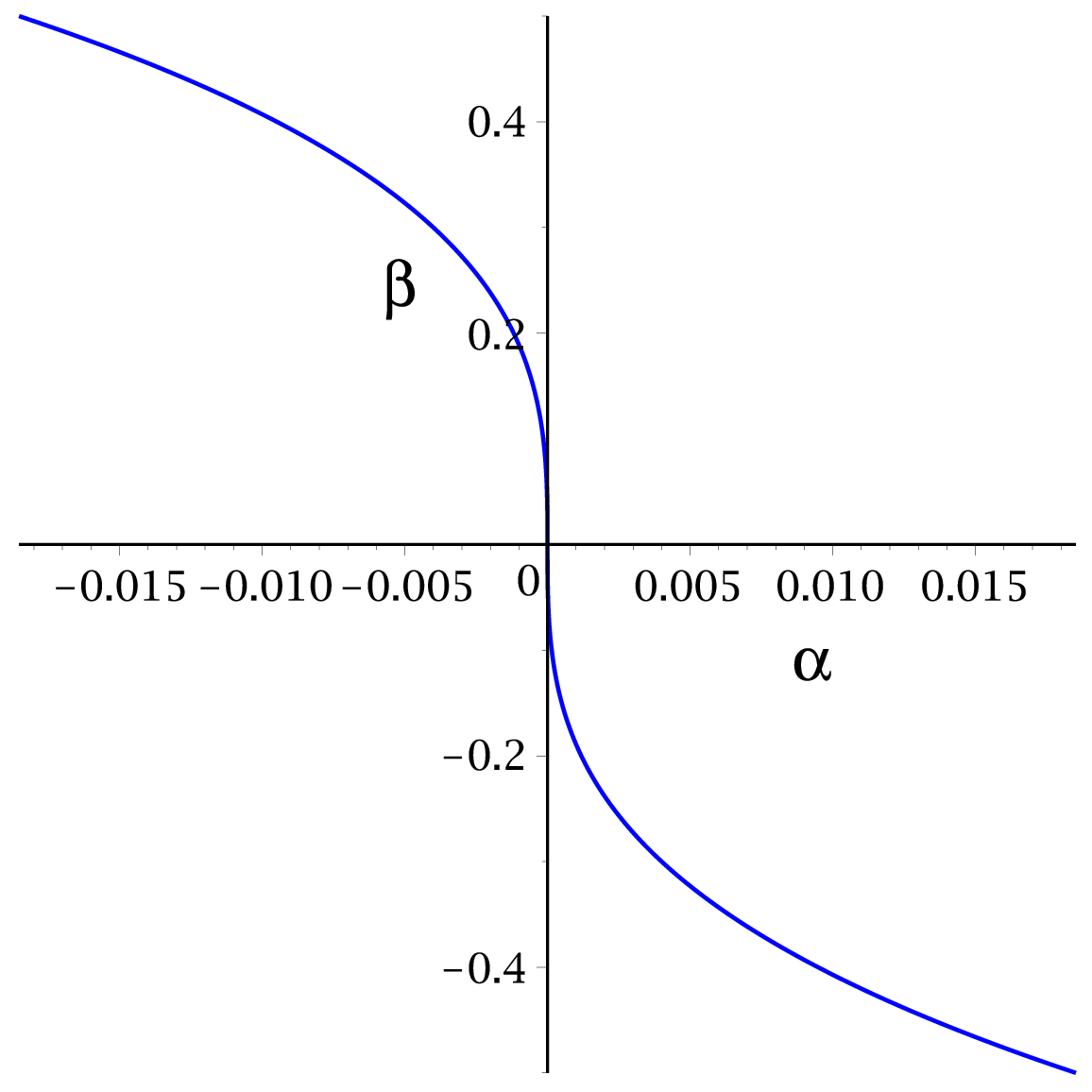}}
\caption{\small The parametric variety over which the input germ \(x^3-\lambda x+\alpha+\beta x^2\) experiences the singularity \(x^2-\lambda\).}
\end{center}
\end{figure}

\item \verb"CheckSingularity"\((x^3-\lambda x+\alpha+\beta x^2, [\alpha, \beta], [x, \lambda]\), \verb"'interval'"='\([\alpha=[-1,1]]\)', \verb"'type'"='\(x^2-\lambda\)', \verb"ParametricVariety") gives

``The bifurcation problem has singularity of \verb"'type'"='\(x^2-\lambda\)' for all values of parameters.''

\item \verb"CheckSingularity"\((x^3-\lambda x+\alpha+\beta x^2, [\alpha, \beta], [x, \lambda]\), \verb"'type'"='\(x^2-\lambda\)', \verb"ParametricVariety", \verb"NonZero") leads to

``Non-zero conditions are \(\lbrace \beta \eta_2 + 1 = 0\rbrace\) where the extra parameters are auxiliary.''

We remark that the equality \(\beta\eta_2+1=0\) for an auxiliary parameter \(\eta_2\) is indeed translated to the inequality \(\beta\neq0\); that is a nonzero condition.
\end{enumerate}
\end{example}

\section{Ring and truncation degree}\label{SecVerify}

In this section we describe how to determine the permissible computational ring and the truncation degree. In fact, each smooth germ with a nonzero-infinite Taylor series expansion must be truncated at certain degree. Further, there are four different options in {\tt Singularity} for computational rings, \ie polynomial, fractional, formal power series and smooth germ rings, that they each can be used by commands in {\tt Singularity} for bifurcation analysis of each singular germ. The command \verb"Verify(g)" derives the following information about the singular germ \(g\) for correct and efficient computations:
\begin{enumerate}
  \item the permissible computational rings for the germ \(g\).
  \item the least permissible truncation degree \(k\) for computations involving the germ \(g\). In other words, the computations modulo degrees of higher than (and equal to) \(k\) does not lead to error.
  \item our recommended computational ring.
\end{enumerate}
We are also interested in the above information for the following purpose:
\begin{itemize}
  \item A list of germs \(G\) in two variables for either division, standard basis computations, multiplication matrix, or intrinsic part of an ideal.
  %\item A parametric germ \(H(x, \lambda, \alpha)\) with \(\alpha\in \mathbb{R}^p\) for either transition set computation or persistent bifurcation diagram classification.\\
\end{itemize}

\begin{tabularx}{\textwidth}{l|X}
  \textbf{Command} & \textbf{Description} \\
\hline
\verb"Verify"(\(g\), \verb"Vars")& derives the permissible computational rings and permissible truncation degree.
\\ \hline
Default upper bound for truncation& a permissible truncation degree is computed as long as it is less than or equal to \(20.\)
\end{tabularx}

\subsection{Options}
\begin{itemize}
  \item \verb"Ideal"; \verb"Persistent"; the command \verb"Verify"(\(G\), \verb"Vars", \verb"Ideal") deals with an ideal \(I\) generated by \(G\), \ie \(I:= \langle G\rangle_\E,\) where \(G\) is a list of germs. It returns the permissible computational ring and a permissible truncation degree when the ideal \(I\) is of finite codimension. Otherwise, it remarks that ``the ideal is of infinite codimension."
%However, the command \verb"Verify"(\(H\), \verb"Vars", \verb"Persistent") determines the least permissible truncation degree \(k\) so that the computations associated with either persistent bifurcation diagram classification or transition sets would be correct.
  \item \verb"Fractional"; \verb"Formal"; \verb"SmoothGerms"; \verb"Polynomial"; the command uses either the rings of fractional germs, formal power series or ring of smooth germs.
  \item Upper bound for truncation degree \(N\); this lets the user to change the default upper bound truncation degree from \(20\) to \(N\).
 \item \verb"'SingularPoint'='[a, b]'"; this option enables the possibility of using the values \((x, \lambda)= (a, b)\)  (when \(x\) stands for the state variable and \(\lambda\) stands for the distinguished parameter in \(g\)) as the base point for the germ \(g.\)
\end{itemize}
\begin{example}
\verb"Verify"(\(x^3-\sin(\lambda), [x, \lambda]\)) gives

The following rings are allowed as the means of computations:

Ring of smooth germs

Ring of formal power series

Ring of fractional germs

The truncation degree must be: 4

\end{example}
\begin{example}
\verb"Verify"(\(x^3-\sin(\lambda), [x, \lambda]\), 2) gives the following warning message:

\begin{center}
``Increase the upper bound for the truncation degree!''
\end{center}
\end{example}
\begin{example}
\verb"Verify"(\([x^4-x\sin(\lambda), x^3\lambda-\lambda\sin(\lambda), 3x^4, 3x^2\lambda], [x, \lambda], \verb"Ideal"\)) gives

The following rings are allowed as means of computations:

Ring of smooth germs

Ring of formal power series

Ring of fractional germs

The truncated degree must be: 4

\end{example}
\begin{example}
Command \verb"Verify"(\(x^3-\sin(\lambda)\), \([x,\lambda]\), \verb"Persistent") gives the least
permissible truncation degree to be 2.
\end{example}
\begin{example}
\verb"Verify"(\(x^3-\sin(\lambda-1), [x, \lambda]\), \verb"'SingularPoint'"='\([0,1]\)') gives

The following rings are allowed as the means of computations:

Ring of smooth germs

Ring of formal power series

Ring of fractional germs

The truncation degree must be: 4
\end{example}
\begin{remark}
We remark that when a point is considered as an input in the computations, our library does not check whether it is singular or not. Thus, the user should use the command {\tt Verify} or {\tt CheckSingularity} in order to derive the associated singular point varieties. 
\end{remark}

\section{Normal form}\label{SecNormalform}
A germ \(f(x, \lambda)\) is called a normal form for the singular germ \(g(x, \lambda)\) when \(f\) has a minimal set of monomial terms in its Taylor expansion among all contact-equivalent germs to \(g\). Therefore, it is easier to analyze the solution set of \(f\) while it has the same qualitative behavior as zeros of \(g\) do.\\

\begin{tabularx}{\textwidth}{l|X}
  \textbf{Command/the default options} & \textbf{Description} \\
\hline
\verb"Normalform"(\(g\),\verb"Vars" )& This function derives a normal form for \(g\).
\\\hline
The computational ring & The default ring is the ring of fractional germs.
\\\hline
Verify/Warning/Suggestion & It automatically verifies if fractional germ ring is sufficient for normal form computation of the input germ \(g.\) Otherwise, it writes a warning note along with a cognitive suggestion for the user.
\\\hline
Truncation degree  & It, by default, detects the maximal degree \(k\) in which \(\mathcal{M}^{k+1}\subseteq \mathscr{P}(g).\) Thus, normal forms are computed modulo degrees higher than or equal to \(k+1.\)
\\\hline
Input germ & It, by default, takes the input germ \(g\) as a polynomial or a smooth germ. It truncates the smooth germs modulo the degree \(k+1.\)
\end{tabularx}

\subsection{Options}
\begin{itemize}
  \item \(k;\) specifies the degree \(k\) so that computations are performed modulo degree \(k+1\). When the input degree \(k\) is too small for the input singular germ \(g\), the computation is not reliable. Thus, an error warning note is returned to inform the user of the situation along with cognitive suggestions to circumvent the problem.
  \item \verb"'SingularPoint'='[a,b]'"; this option enables the possibility of using a base point different from zero values for both the state variable
and the distinguished parameter.
  \item \verb"Fractional", \verb"Formal", \verb"SmoothGerms"; \verb"Polynomial"; the command uses either the rings of fractional germs, formal power series or ring of smooth germs. When it is necessary, warning/cognitive suggestions are given accordingly.
  \item \verb"list"; this generates a list of possible normal forms for the germ \(g\). Different normal forms may only occur due to possible alternative eliminations in intermediate order terms.
\end{itemize}
\begin{example}

\verb"Normalform"(\(x^3-\sin(\lambda), [x, \lambda]\), 10, \verb"SmoothGerms") generates
\bes
x^3-\lambda.
\ees while  \verb"Normalform"(\(1-\frac{1}{1+x^4-\lambda^2}, [x, \lambda]\), 10, \verb"Formal") gives rise to
\bas
x^4-\lambda^2.
\eas Using \verb"Normalform"(\(x^5+x^3\lambda+\sin(\lambda^2), [x, \lambda]\), 10, \verb"Polynomial") gives the following suggestion and warning note.

Warning: The polynomial germ ring is not suitable for normal form computations.

Suggestion: Use the command \verb"Verify" to find the appropriate computational ring.

The following output might be wrong.

The germ is an infinite codimensional germ.

In fact the above statement is wrong since the high order term ideal contains
\bes
\mathcal{M}^6+\mathcal{M}^4\langle \lambda\rangle+\mathcal{M}\langle \lambda^2\rangle.
\ees

\noindent Now the command \verb"Verify"(\(x^5+x^3\lambda+\sin(\lambda^2), [x, \lambda])\) gives rise to

Fractional germ ring; Formal power series ring; Smooth germ ring.

\noindent Thus, we use \verb"Normalform"(\(x^5+x^3\lambda+\sin(\lambda^2)\), \([x, \lambda]\), 10, \verb"Fractional") to obtain

\bes
x^5+x^3\lambda+\lambda^2.
\ees
\end{example}

\section{Universal unfolding}\label{SecUniversalUnfolding}

Generally in dealing with singular problems, extra complications are experienced in the laboratory data than what are predicted by the modeling theoretical analysis. The problem here is due to {\it modeling imperfections}; natural phenomena can not be perfectly modeled by a mathematical model. In fact one usually neglects the impact of many factors like friction, pressure, and/or temperature, etc., to get a manageable mathematical model. Otherwise one will end up with a mathematical modeling problem with too many or infinite number of parameters. The imperfections around singular points may cause dramatic qualitative changes in the solution set of the model. Universal unfolding gives us a natural way to circumvent the problem of imperfections.

\begin{defn}\label{Def2.3}
A parametric germ \(G(x, \lambda, \alpha)\) for \(\alpha\in \mathbb{R}^m\) is called an {\em unfolding} for \(g(x, \lambda)\) when
\begin{equation*}
G(x, \lambda, 0)=g(x, \lambda).
\end{equation*} An unfolding \(G(x, \lambda, \alpha)\) for \(g\) is called a {\em versal unfolding} when for each unfolding \(H(x, \lambda, \beta)\) of \(g(x, \lambda)\) there is a smooth germ \(\alpha(\beta)\) so that \(H\) is contact-equivalent to \(G(x, \lambda, \alpha(\beta)).\) Roughly speaking, a versal unfolding is a parametric germ that contains a contact-equivalent copy of all small perturbations of \(g(x, \lambda)\). A versal unfolding with insignificant parameters is not suitable for the bifurcation analysis. So, we are interested in a versal unfolding that has a minimum possible number of parameters, that is called {\em universal unfolding}. In other words, universal unfolding has the minimum possible number of parameters so that they accommodate all possible qualitative types that small perturbations of \(g(x, \lambda)\) may experience.
\end{defn}

\begin{tabularx}{\textwidth}{l|X}
  \textbf{Command/the default options} & \textbf{Description} \\
\hline
 \verb"UniversalUnfolding"(\(g\),\verb"Vars" )& This function computes a universal unfolding for \(g\).
\\\hline
The computational ring & By default, {\tt Singularity } uses the ring of fractional germs.
\\\hline
Verify/Warning/Suggestion & This automatically derives the least sufficient degree for truncations and also verifies if fractional germ ring is sufficient for computation. Otherwise, it writes a warning note along with guidance on the suitable rings for computations and hints at other possible capabilities of {\tt Singularity.}
\\\hline
Degree  & It, by default, detects the maximal degree \(k\) in which terms of degree higher than or equal to \(k+1\) can be ignored. Thus, the computations are performed modulo degree \(k+1.\)
\\\hline
Input germ & The default input germ \(g\) is a polynomial or a smooth germ. For an input smooth germ, the default procedure \verb"UniversalUnfolding" truncates the smooth germs modulo \(k+1,\) \ie modulo degrees higher than and equal to \(k+1.\)
\end{tabularx}

\subsection{Options}
\begin{itemize}
  \item \verb"normalform"; A universal unfolding for normal form of \(g\) is derived by this option.
  \item \verb"list"; this function provides the list of possible universal unfoldings for \(g\).
  \item \(k\); the degree \(k\) determines the truncation degree so that all computations are performed modulo \(k+1.\) For low degrees of \(k,\) it may
  derive wrong results. Thus, it gives a warning error and a suggestion for the user when \(k\) must be a larger number for correct result.
  \item \verb"Fractional"; \verb"Formal"; \verb"SmoothGerms"; \verb"Polynomial"; this determines the computational ring. The command
  \verb"UniversalUnfolding" gives a warning note when the user's choice of computational ring is not suitable for computations involving
  the input germ \(g\) and writes a suggestion to circumvent the problem.
\item \verb"'SingularPoint'"='\([a,b]\)'; this provides numerical values \((a, b)\) for the state variable and the
distinguished parameter, respectively.

\end{itemize}

\begin{example} \verb"UniversalUnfolding"\((x^4+4x^3-\lambda x, [x, \lambda], \verb"normalform", \verb"list")\) gives rise to
\bas
& x^3-x\lambda+ \alpha_1+\alpha_2 \lambda &\\
& x^3-x\lambda+ \alpha_1+\alpha_2 x^2.&
\eas \verb"UniversalUnfolding"\((6x-6\sin(x)-\lambda x, [x, \lambda], 6, \verb"normalform", \verb"list", \verb"Formal")\) leads to
\bas
& x^3-\lambda x+\alpha_{2}x^2+\alpha_{1} &\\
& x^3-\lambda x+\alpha_{2}\lambda+\alpha_{1}.&
\eas
Now consider \(g(x, \lambda):=x^6+x^4\lambda+\lambda^2 .\)

\verb"UniversalUnfolding"\((g(x, \lambda), [x, \lambda], 6, \verb"normalform", \verb"list", \verb"Polynomial")\) gives the following warning error and suggestion:

\noindent \emph{Warning: The ring of polynomial germs is not suitable for normal form computations of \(g.\) }

\noindent \emph{Suggestion: The permissible computational ring options are} \verb"Fractional", \verb"SmoothGerms" \emph{and} \verb"Formal".
\end{example}
\subsection{Universal unfolding parameters}

Consider a given parametric germ \(G(x, \lambda, \alpha)\) in defenition \ref{Def2.3} when \(m\) is greater than the codimension of \(g.\) Therefore, it is helpful to find out all possible list of parameters which they can play the role of universal unfolding parameters. This is handled in {\tt Singularity} through the following function.

\begin{tabularx}{\textwidth}{l|X}
  \textbf{Command} & \textbf{Description} \\
\hline
\verb"UniversalUnfoldingPars"(\(G\), \verb"Pars", \verb"Vars") & This function finds the list of
parameters which they can play the role of universal unfolding parameters.
\end{tabularx}
\subsubsection{Options}
\begin{itemize}
\item \verb"'SingularPoint'"='\([a,b]\)'; this provides numerical values \((a, b)\) for the state variable and the
distinguished parameter, respectively.
\end{itemize}
\begin{example}
The command
\verb"UniversalUnfoldingPars" \((x^3-\lambda x+\sin(\alpha_1)+\exp(\alpha_2\lambda+\alpha_3)+\alpha_4 x^2+\alpha_5\cos(x\lambda+\lambda)-1, [\alpha_1, \alpha_2, \alpha_3, \alpha_4, \alpha_5],[x, \lambda])\) gives

$$[\alpha_1, \alpha_2], [\alpha_1, \alpha_4], [\alpha_2, \alpha_3], [\alpha_2, \alpha_5], [\alpha_3, \alpha_4], [\alpha_4, \alpha_5].$$

This provides six alternative universal unfoldings for the pitchfork problem \(x^3-\lambda x\).
\end{example}
\section{Recognition problem }\label{SecRecognitionProblem}

We describe the command \verb"RecognitionProblem" on how it answers the recognition problem, that is, what kind of germs have the same normal form or universal unfolding for a given germ \(g\)?
\subsection{Normal form}
\begin{itemize}
  \item \textbf{Low order terms}. Low order terms refer to the monomials in
\(\s(g)^{\perp}\) which do not appear in any contact-equivalent copy of \(g\).
  \item \textbf{High order terms}. The ideal \(\p(g)\) represents the space of negligible terms that are called high order terms. These terms are eliminated in normal form of \(g.\)
  \item \textbf{Intermediate order terms}. A monomial term are called an intermediate order term when it is neither low order nor high order term. Intermediate order terms may or may not be simplified in normal form computation of smooth germs.
\end{itemize}
The answer for the recognition problem for normal form of a germ \(g\) is a list of zero and nonzero conditions for certain derivatives of a hypothetical germ \(f.\) When these zero and nonzero conditions are satisfied for a given germ \(f,\) the germ \(f\) and \(g\) are contact-equivalent. Each germ with a minimal list of monomial terms in its Taylor expansion constitutes a normal form for \(g.\)
\subsection{Universal unfolding}\label{SecCheckUniversal}

Consider a parametric germ \(G(x, \lambda, \alpha)\) and a germ \(g(x, \lambda).\) Then, \(G\) is usually a universal unfolding for \(g\) when \(G(x, \lambda, 0)= g(x, \lambda)\) and certain matrix associated with \(G\) has a nonzero determinant. Thus, the answer of the recognition problem for universal unfolding is actually a matrix whose components are derivatives of a hypothetical parametric germ \(G\) satisfying \(G(x, \lambda, 0)= g(x, \lambda)\).\\

\begin{tabularx}{\textwidth}{l|X}
  \textbf{Command/ default} & \textbf{Description} \\
\hline
 \verb"RecognitionProblem"(\(g\), \verb"Vars") & returns a list of zero and nonzero conditions on certain derivatives of a hypothetical germ \(f\). A given germ \(f\) is contact-equivalent to \(g\) when those conditions are satisfied. \\
 \hline
Computational ring & The default is fractional germ ring. The lack of warning notes is a confirmation that the fractional germ ring is suitable for computation.
\\
  \hline
Truncation degree & It automatically computes an optimal truncation degree \(k\) and performs the remaining computations modulo degrees of higher than (but not  equal to) \(k\). \\
\hline
Verification/Warning& A warning note of possible errors is given when the computational ring is not suitable for the germ \(g.\) The truncation degree is also checked and if it is not sufficiently large enough, a warning note is given. Warning notes are accompanied with cognitive suggestions to circumvent the problem.
\end{tabularx}

\subsubsection{Options}
\begin{itemize}
\item \(k\); this number represents the truncation degree.
\item Computational ring: \verb"Fractional", \verb"Formal", \verb"SmoothGerms"; \verb"Polynomial"; the command accordingly uses either the rings of
fractional germs, formal power series or smooth germs.
\item \verb"UniversalUnfolding"; it returns a matrix. The matrix components consists of
certain derivatives of a hypothetical parametric germ \(G\). Then, a parametric germ \(G(x, \lambda, \alpha)\) is a universal unfolding for \(g(x, \lambda)\) when \(G(x, \lambda, 0)= g(x, \lambda)\) and the associated matrix has a nonzero determinant.
\item \verb"subs"; it returns the same output as the previous option except that the derivatives of \(g\) are replaced with their numeric values.
\item \verb"'SingularPoint'"='\([a,b]\)'; this option refers to the base point given as \(a\) for the state variable and \(b\) for the
distinguished parameter.
\end{itemize}

\begin{example}
\verb"RecognitionProblem"(\(x^3+\sin(\lambda), [x, \lambda]\), 6, \verb"Formal") gives rise to
\begin{center}
 "nonzero condition=", \([\frac{\partial}{\partial\lambda}f\neq 0, \frac{\partial^3}{\partial x^3}f\neq 0]\)\\
"zero condition=", \([f=0, \frac{\partial}{\partial x}f=0,\frac{\partial^2}{\partial x^2}f =0]\).
\end{center}

\noindent \verb"RecognitionProblem"(\(x^3+\exp(\lambda^2)-1, [x, \lambda]\), 6, \verb"UniversalUnfolding", \verb"SmoothGerms") gives rise to
\bes
\det\left( {\begin{array}{ccccc}
0 & 0 & 0 & g_{x,x,x}(0) & g_{x,x,\lambda}(0)\\
0 & g_{\lambda, \lambda}(0) & 0 & g_{x,x,\lambda}(0) & g_{x, \lambda, \lambda}(0)\\
 G_{\alpha_1}(0) & G_{\lambda, \alpha_1}(0) & G_{x, \alpha_1}(0) & G_{x, x, \alpha_1}(0) & G_{x, \lambda, \alpha_1}(0)\\
 G_{\alpha_2}(0) & G_{\lambda, \alpha_2}(0) & G_{x, \alpha_2}(0) & G_{x, x, \alpha_2}(0) & G_{x, \lambda, \alpha_2}(0)\\
  G_{\alpha_3}(0) & G_{\lambda, \alpha_3}(0) & G_{x, \alpha_3}(0) & G_{x, x, \alpha_3}(0) & G_{x, \lambda, \alpha_3}(0)\\
\end{array}}
\right)\neq 0.
\ees\\

Now the command \verb"RecognitionProblem"(\(x^3+\exp(\lambda^2)-1, [x, \lambda]\), 6, \verb"UniversalUnfolding", \verb"Fractional",\verb"subs") returns
\bes
\det\left( {\begin{array}{ccccc}
0 & 0 & 0 & 6 & 0\\
0 & 2 & 0 & 0 & 0\\
 G_{\alpha_1}(0) & G_{\lambda, \alpha_1}(0) & G_{x, \alpha_1}(0) & G_{x, x, \alpha_1}(0) & G_{x, \lambda, \alpha_1}(0)\\
 G_{\alpha_2}(0) & G_{\lambda, \alpha_2}(0) & G_{x, \alpha_2}(0) & G_{x, x, \alpha_2}(0) & G_{x, \lambda, \alpha_2}(0)\\
  G_{\alpha_3}(0) & G_{\lambda, \alpha_3}(0) & G_{x, \alpha_3}(0) & G_{x, x, \alpha_3}(0) & G_{x, \lambda, \alpha_3}(0)\\
\end{array}}
\right)\neq 0.
\ees\\
\end{example}

%\subsection{Check Universal unfolding}

\noindent Let \(G(x,\lambda,\alpha)\) be a parametric germ where \(\alpha\in\mathbb{R}^{k}\).\\

\begin{tabularx}{\textwidth}{l|X}
  \textbf{Command} & \textbf{Description} \\
\hline
\verb"CheckUniversal"(\(G\), \verb"Pars", \verb"Vars") & This function checks if a parametric germ \(G\) is a universal unfolding for
\(G(x,\lambda,0)\)
\end{tabularx}
\begin{itemize}
\item \(k\); this number indicates the truncation degree.
\item \verb"'SingularPoint'='[a,b]'"; this means that the base point of singularity is given as \(a\) for the state variable and \(b\) for the
distinguished parameter.
\end{itemize}
\begin{example}
\verb"CheckUniversal"(\(x^5-\lambda+\alpha_{1}x+\alpha_{2}x^2+\alpha_{3}x^3, [\alpha_{1}, \alpha_{2}, \alpha_{3}], [x, \lambda]\)) gives

\begin{center}
"Yes"
\end{center}
\end{example}

\section{Transformations }\label{SecTransformation}

For each two contact-equivalent germs \(f\) and \(g,\) there are diffeomorphic germs \((X(x, \lambda), \Lambda(\lambda))\) and smooth germ \(S(x, \lambda)>0\) such that \(X_{x}(x, \lambda),\) \(\Lambda^{\prime}(\lambda)>0\) and \(f(x, \lambda)= S(x, \lambda)g(X(x, \lambda), \Lambda(\lambda))\).\\

\begin{tabularx}{\textwidth}{l|X}
  \textbf{Command/option} & \textbf{Description} \\
\hline
\verb"Transformation"(\(g\), \verb"Vars") & This function computes the smooth germs \(X, \Lambda, S\) transforming the germ \(g\) into its normal form modulo degree \(k,\) where terms of degree higher than or equal to \(k\) are high order terms. \\
\hline
Transformation(\(g, f\)) & This function computes suitable smooth maps \(X, \Lambda, S\) for transforming the germ \(g\) into \(f\) modulo high order terms. \end{tabularx}
\subsection{Options}
\begin{itemize}
\item \(k;\) this number specifies a degree \(k\) so that computations are performed modulo degrees of higher or equal to \(k\). When \(k\) is less than the
degrees of high order terms a warning note is given.
\end{itemize}

\begin{example}
\verb"Transformation"(\(x^3+\sin(\lambda)+\exp(x^5)-1, x^3+\lambda, [x, \lambda], 4\)) gives rise to
\bas
X&=&x+\lambda+x \lambda+\lambda^2, \qquad \Lambda(\lambda)= \lambda, \\
S&=&1-3x^2-3x \lambda-\frac{5}{6}\lambda^2-3x^3-9 \lambda x^2-9x \lambda^2-3 \lambda^3.
\eas
\end{example}

\section{Bifurcation diagrams }

Bifurcation diagram analysis of a parametric system is performed by the notion of {\it persistent} and {\it non-persistent} bifurcation diagrams.
Bifurcation diagram of \eqref{eq:1} is defined by
\begin{eqnarray*}
\lbrace (x, \lambda)\mid g(x, \lambda, \alpha)=0 \rbrace.
\end{eqnarray*} A bifurcation diagram is called {\it persistent} when the bifurcation diagrams subjected to small perturbations in parameter space remain self contact-equivalent.

\subsection{Persistent bifurcation diagram classification and transition set }

The classification of persistent bifurcation diagrams are performed by the notion of transition sets. In fact, a subset of parameter space is called {\it transition set} when the associated bifurcation diagrams are non-persistent. {\it Transition set} is denoted by \(\Sigma\) and is usually a hypersurface of codimension one for germs of finite codimension. Then, one choice from each connected components of the complement of the transition set \(\Sigma\) makes a complete persistent bifurcation diagram classification of a given parametric germ. This provides a comprehensive insight into the persistent zero solutions of a parametric germ.

\subsubsection{Transition set }\label{SecTransitionSet}
The parameters associated with non-persistent bifurcation diagrams are split into three categories: {\it bifurcation}, {\it hysteresis}, and {\it double limit point}. These are defined and denoted by
\begin{eqnarray*}
\mathscr{B}&=& \{\alpha\in\mathbb{R}^{p}\mid G=G_{x}=G_{\lambda}=0 \hbox{ at }  (x, \lambda, \alpha) \hbox { for some }(x, \lambda)\in \mathbb{R}\times \mathbb{R}\},\\
 \mathscr{H}&=& \{\alpha\in\mathbb{R}^{p}\mid G=G_{x}=G_{xx}=0 \hbox{ at } (x, \lambda, \alpha) \hbox{ for some }(x, \lambda)\in \mathbb{R}\times \mathbb{R}\},\\
 \mathscr{D}&=& \{\alpha\in\mathbb{R}^{p}\mid\exists(x_{1}, x_{2}, \lambda)\in \mathbb{R}\times \mathbb{R}\times \mathbb{R}\hbox{ so that } G =G_{x}=0 \hbox{ at } (x_{1},\lambda,\alpha) \hbox{ and }\\ &&\; (x_{1},\lambda,\alpha) \}.
\end{eqnarray*}
The transition set \(\Sigma\) is now given by \(\Sigma:=\mathscr{B}\cup\mathscr{H}\cup\mathscr{D}\).
Suppose that \(H\) is a singular parametric germ. \\

\begin{tabularx}{\textwidth}{l|X}
  \textbf{Command/the default options} & \textbf{Description} \\
\hline
  \verb"TransitionSet"(\(H\), \verb"Pars", \verb"Vars")& This function estimates the transition set in terms of parameters of \(H.\) The default is to eliminate \(x\) and \(\lambda\) variables from the equations given by \(\mathscr{B}, \mathscr{H}, \mathscr{D}\).
  \\\hline
Truncation degree & For non-polynomial input germs, by default, it automatically computes a suitable truncation degree \(k\) and truncates the input germ at degree \(k\), \ie preserving degrees of less than and equal to \(k\).
\end{tabularx}
\subsubsection{Options}
\begin{itemize}
  \item \verb"Pars":=\([\alpha_1, \alpha_2, \ldots, \alpha_p]\); this hints to derive the transition set in terms of these variables while the attempts are best made to eliminate the rest of variables from the equations (as many as possible).
  \item \verb"plot"; this function plots/animates transition set in parameter space.
  \item \(\alpha_{i}\); When codimension is more than or equal to three 3, some parameters \(\alpha_3, \alpha_4, \ldots\) will be taken as fixed by default. This option refines the \verb"plot"/\verb"animate" option by allowing to change the fixed parameters to \(\alpha_{i}. \)
  \item \(k\); determines the truncation degree. The user may use the command \verb"Verify" to find an appropriate degree \(k.\)
\item \verb"'SingularPoint'='[a,b]'"; the computations are performed for the base point given as \(a\) for the state variable and \(b\) for its distinguished parameter.
\end{itemize}
\begin{example}
Here, we bring two examples from \cite[Page 206]{GolubitskySchaefferBook}.

\verb"TransitionSet"(\(x^4-\lambda x+\alpha_1+\alpha_2 \lambda+\alpha_3 x^2, [\alpha_1, \alpha_2, \alpha_3], [x, \lambda]\)) gives rise to
\bas
\mathscr{B}&:=&\{(\alpha_1, \alpha_2, \alpha_3)\,|\,\alpha_2^4+\alpha_2^2\alpha_3+\alpha_1=0\},\\
\mathscr{H}&:=&\{(\alpha_1, \alpha_2, \alpha_3)\,|\,128\alpha_2^2\alpha_3^3+3\alpha_3^4+72\alpha_1\alpha_3^2+432\alpha_1^2=0\},\\
\mathscr{D}&:=&\{(\alpha_1, \alpha_2, \alpha_3)\,|\,-\alpha_3^2+4\alpha_1=0,\,\alpha_3\leq 0\}.
\eas

\verb"TransitionSet"(\(x^5-\lambda+\alpha_1 x+\alpha_2 x^2+\alpha_3 x^3, [ \alpha_1, \alpha_2, \alpha_3], [x, \lambda]\)) derives
\bas
\mathscr{B}&:=&\emptyset,\\
\mathscr{H}&:=&\{(\alpha_1, \alpha_2, \alpha_3)\,|\,-81\alpha_1 \alpha_3^4+27\alpha_2^2\alpha_3^3+360\alpha_1^2\alpha_3^2\\
&&-540\alpha_1\alpha_2^2\alpha_3+135
\alpha_2^4-400\alpha_1^3=0\},\\
\mathscr{D}&:=&\{(\alpha_1, \alpha_2, \alpha_3)\,|\,-16\alpha_3^6+224\alpha_1\alpha_3^4-88\alpha_2^2\alpha_3^3\\&&-1040\alpha_1^2\alpha_3^2+360
\alpha_1\alpha_2^2\alpha_3+
135\alpha_2^4+1600\alpha_1^3=0\}.
\eas
\end{example}

\begin{example}
\verb"TransitionSet"(\(x^3+\sin(\lambda x)+\alpha_1+\alpha_2 x^2, [\alpha_1, \alpha_2], [x, \lambda], 5\), \verb"plot") generates Figure \ref{1}.
\begin{figure}[h]
\begin{center}
\subfigure[\label{1}]{\includegraphics[width=.40\columnwidth,height=.4\columnwidth]{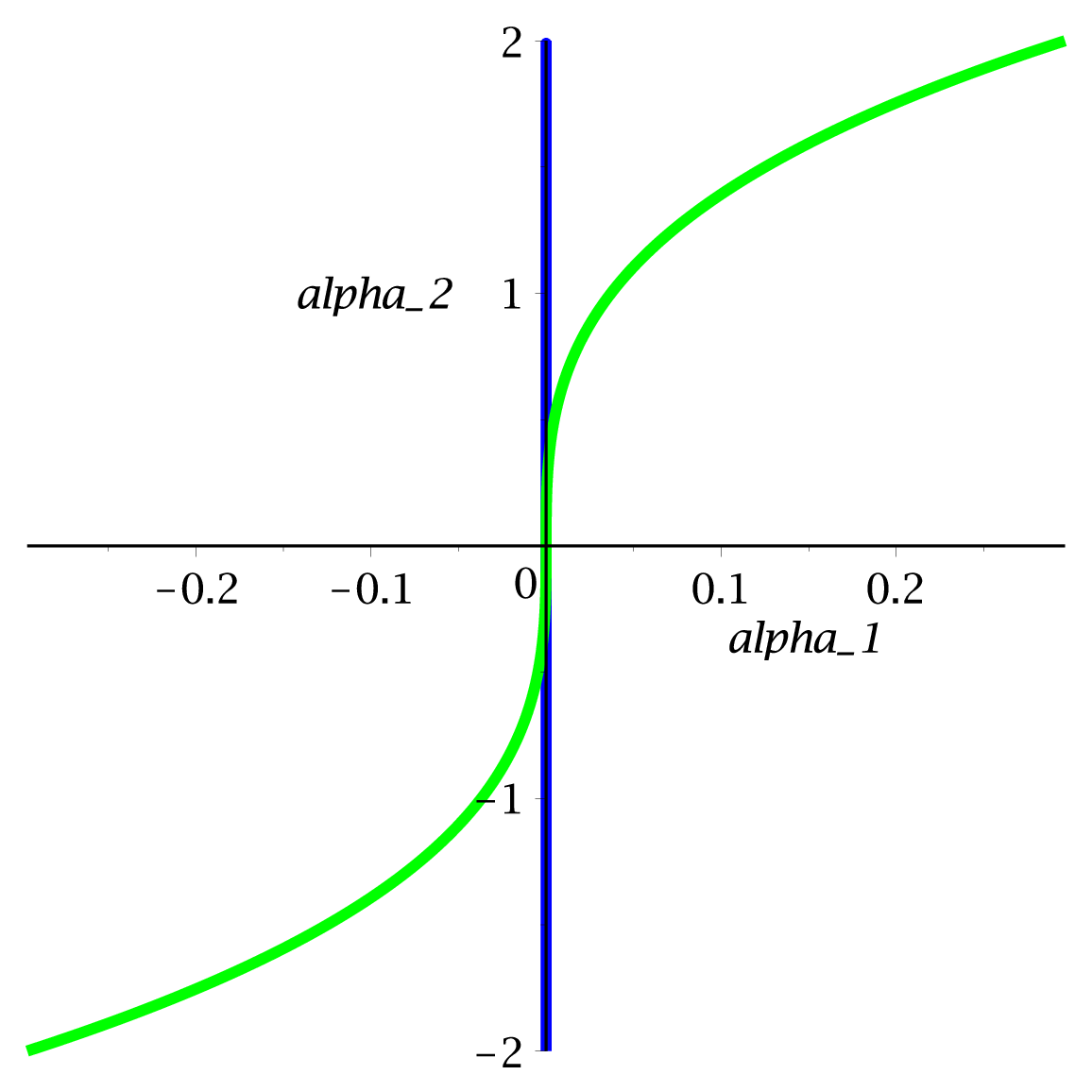}}
\subfigure[\label{2}]{\includegraphics[width=.40\columnwidth,height=.4\columnwidth]{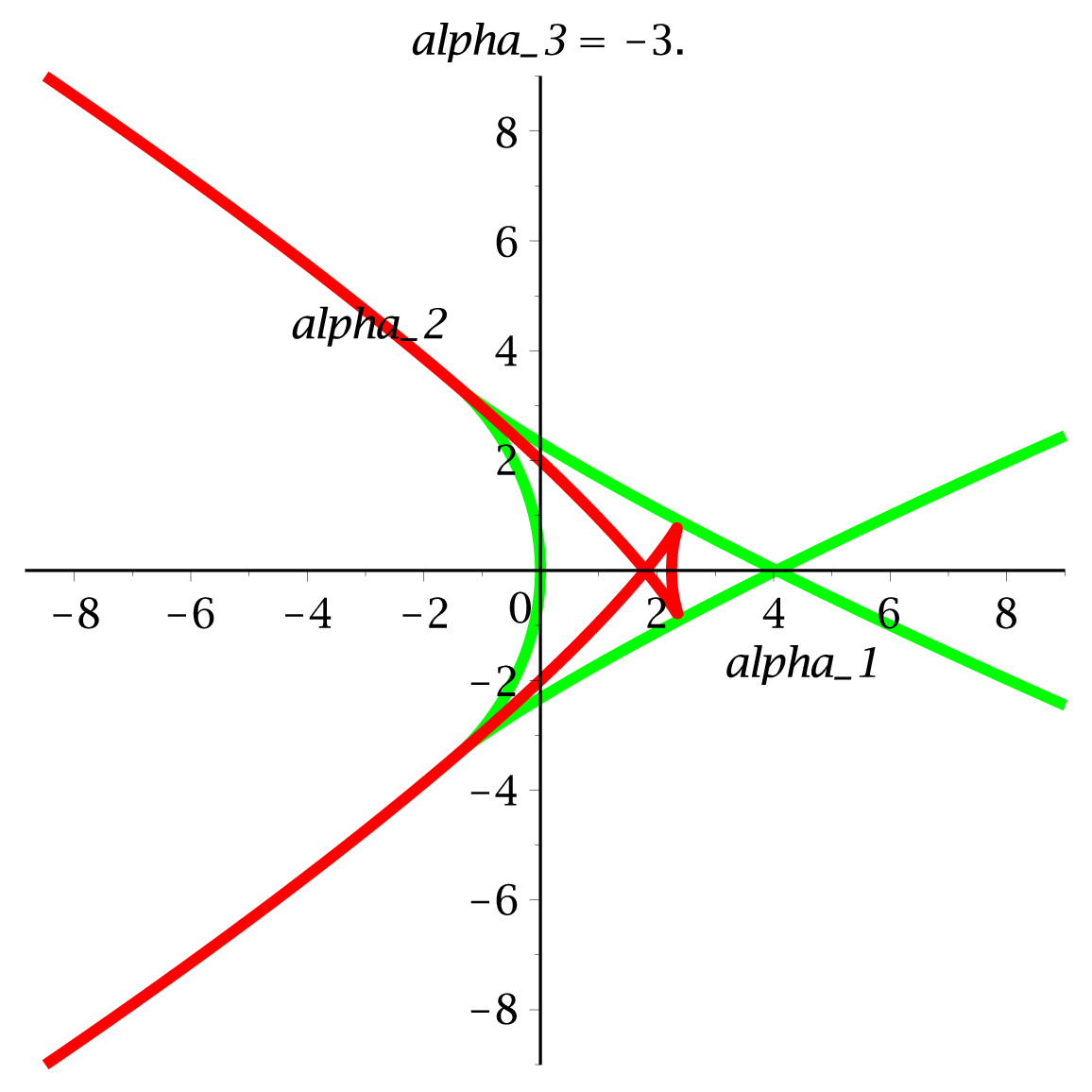}}
\caption{\small Blue color is for bifurcation \(\mathscr{B}\), green stands for hysteresis \(\mathscr{H}\) while the red color depicts transition set associated with double limit point \(\mathscr{D}\).}
\end{center}
\end{figure}

\verb"TransitionSet"(\(x^5-\lambda+\alpha_1 x+\alpha_2 x^2+\alpha_3 x^3, [\alpha_1, \alpha_2, \alpha_3], [x, \lambda], \verb"plot"\)) creates an animation ending at Figure \ref{2}.
\end{example}

\subsubsection{Persistent bifurcation diagrams} \label{SecPersistentDiagram}
The command \(\verb"PersistentDiagram"\) follows the following table.\\

\begin{tabularx}{\textwidth}{l|X}
  \textbf{Command/the default options} & \textbf{Description} \\
\hline
\verb"PersistentDiagram"(\(H\), \verb"Vars") & this function plots/animates bifurcation diagrams in \(x-\lambda\) plane by passing through the parameter space. The default path is a circle around the singular point and it may include at most two parameters.
\end{tabularx}

\subsubsection{Options}
\begin{itemize}
  \item \(\alpha_{i}\); for parameter space of dimension more than two, it chooses three values for \(\alpha_i\), \ie a negative, zero and a positive value for \(\alpha_{i}\). Then, it plots/animates bifurcation diagrams in \(x-\lambda\) plane by passing (circular path by default) through parameter space for fixed parameter \(\alpha_i\).
  \item \(f(\zeta), g(\zeta)\); this function animates bifurcation diagrams in \(x-\lambda\) plane by passing through the given path (\(f(\zeta), g(\zeta)\)).
  \item \verb"ShortList"; \verb"IntermediateList"; \verb"CompleteList"; either of these generates a list of points associated persistent bifurcation diagrams.
  \item \verb"plot"; this option plots the persistent bifurcation diagrams associated with parameter points output of the previous option, \ie \verb"ShortList"; \verb"IntermediateList"; \verb"CompleteList".
  \item \(k\); determines the truncation degree of the germ \(H.\)
  \item \verb"'values'"='\([\alpha_1=a_1, \ldots, \alpha_m=a_m]\)'; plots persistent diagrams where parameters \(\alpha_i\) receive the values \(a_i\) for \(i=1, \ldots, m\).
  \item \verb"'IntervalPlot'"='\([x=x_0..x_1, \lambda=\lambda_0..\lambda_1]\)'; plots persistent diagrams over the vertical and horizontal axes \(x\) and \(\lambda\) with the ranges
  \(x_0..x_1\) and \(\lambda_0..\lambda_1\), respectively.
\end{itemize}

\begin{example}
Now we present how the command \verb"PersistentDiagram" works.

\verb"PersistentDiagram"(\(x^5-\lambda+\alpha_1 x+\alpha_2 x^2+\alpha_3 x^3\), \([x, \lambda]\) , \verb"plot",  \verb"CompleteList") generates a list from which the list of inequivalent bifurcation diagrams are chosen in Figure \ref{Fig2}.

\begin{figure}[h]
\begin{center}
\subfigure[\label{a}]{\includegraphics[width=.26\columnwidth,height=.2\columnwidth]{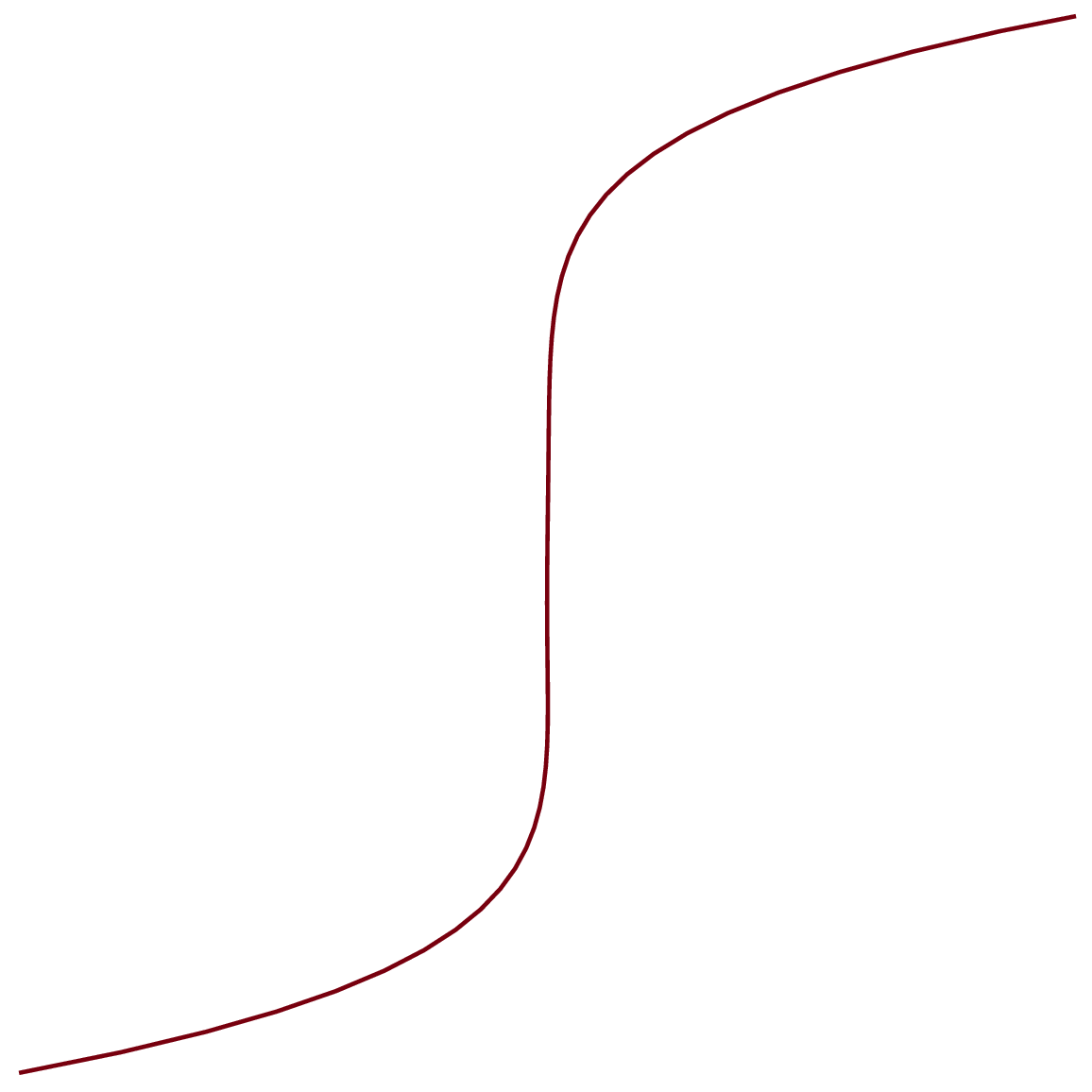}}
\subfigure[\label{b}]{\includegraphics[width=.26\columnwidth,height=.2\columnwidth]{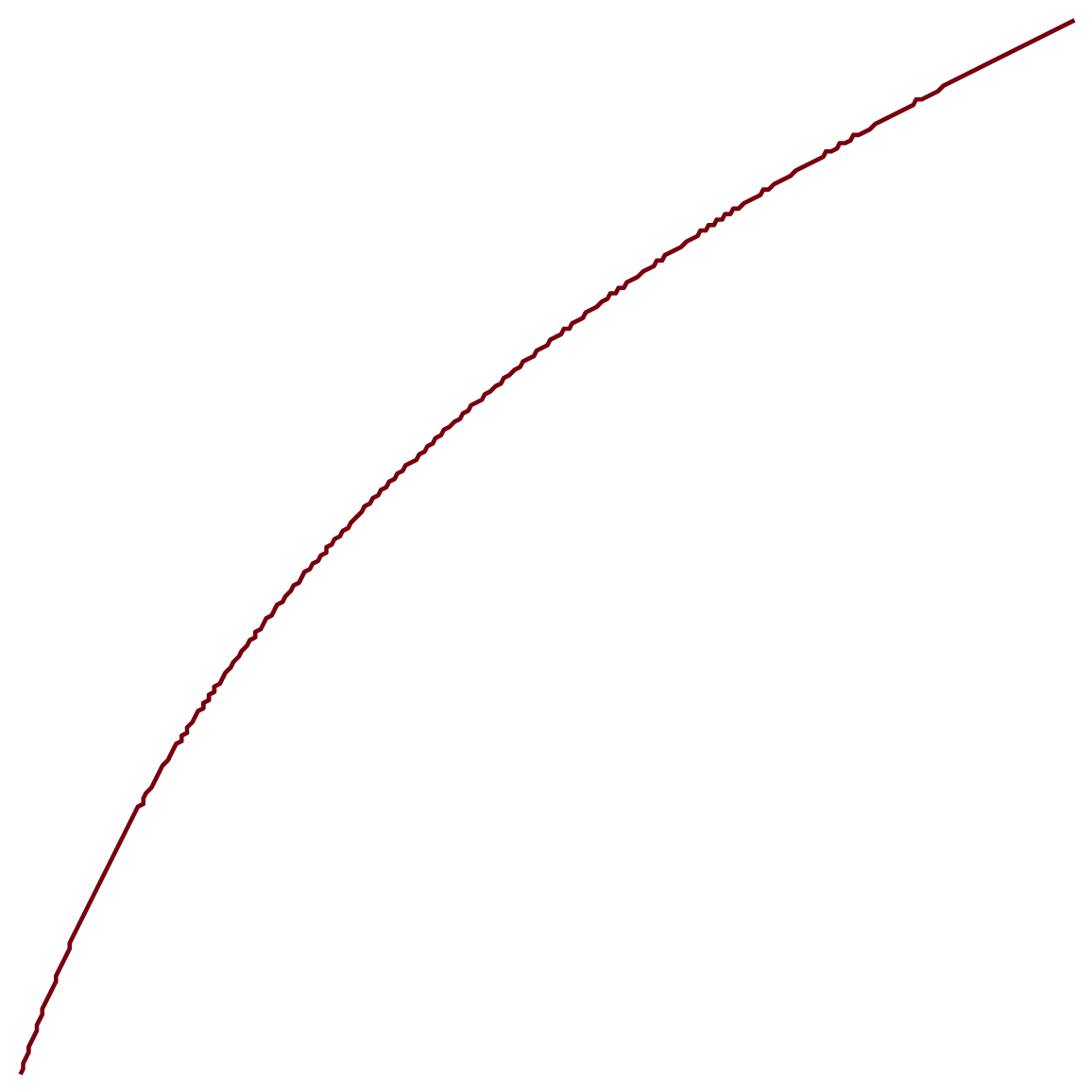}}
\subfigure[\label{c}]{\includegraphics[width=.26\columnwidth,height=.2\columnwidth]{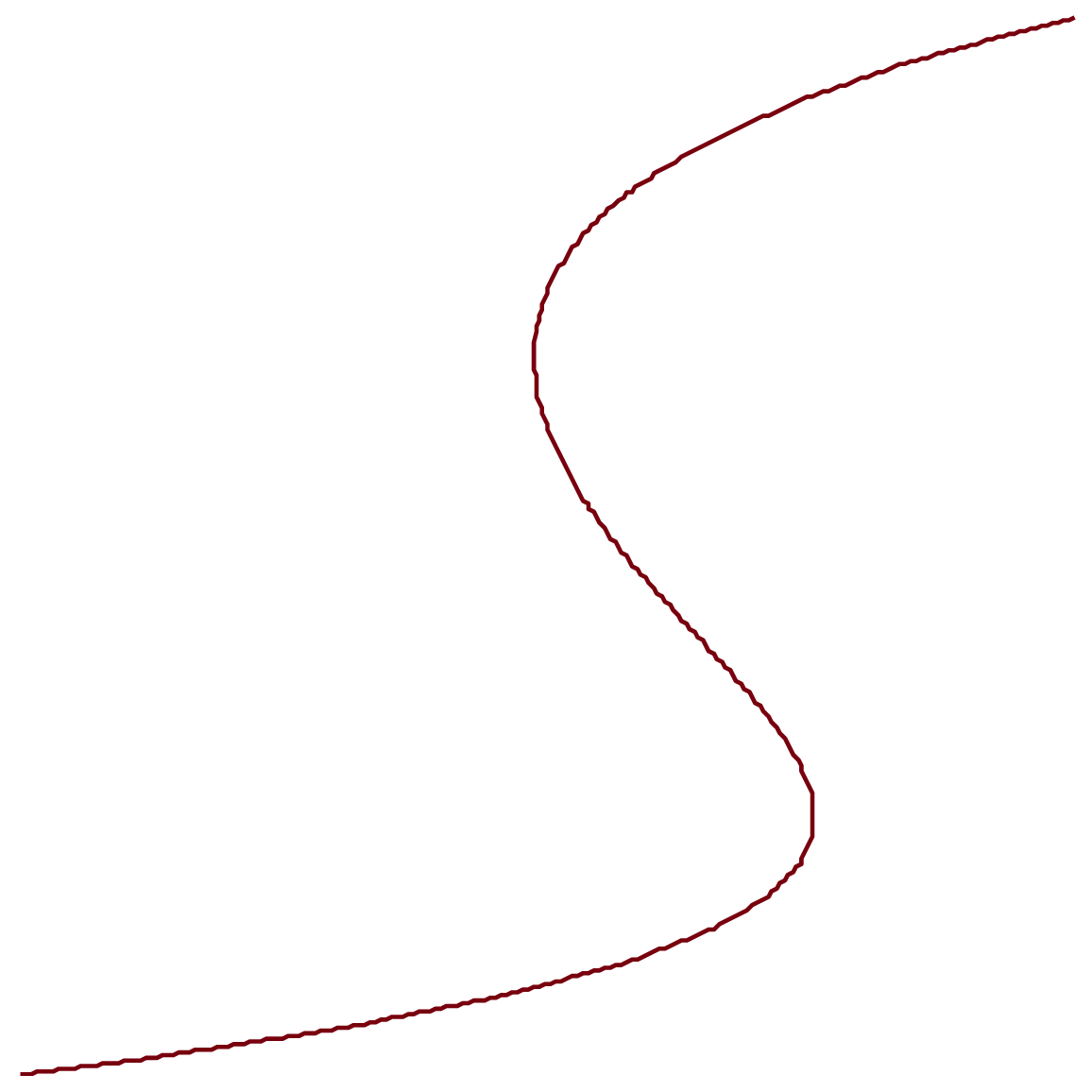}}
\subfigure[\label{d}]{\includegraphics[width=.26\columnwidth,height=.2\columnwidth]{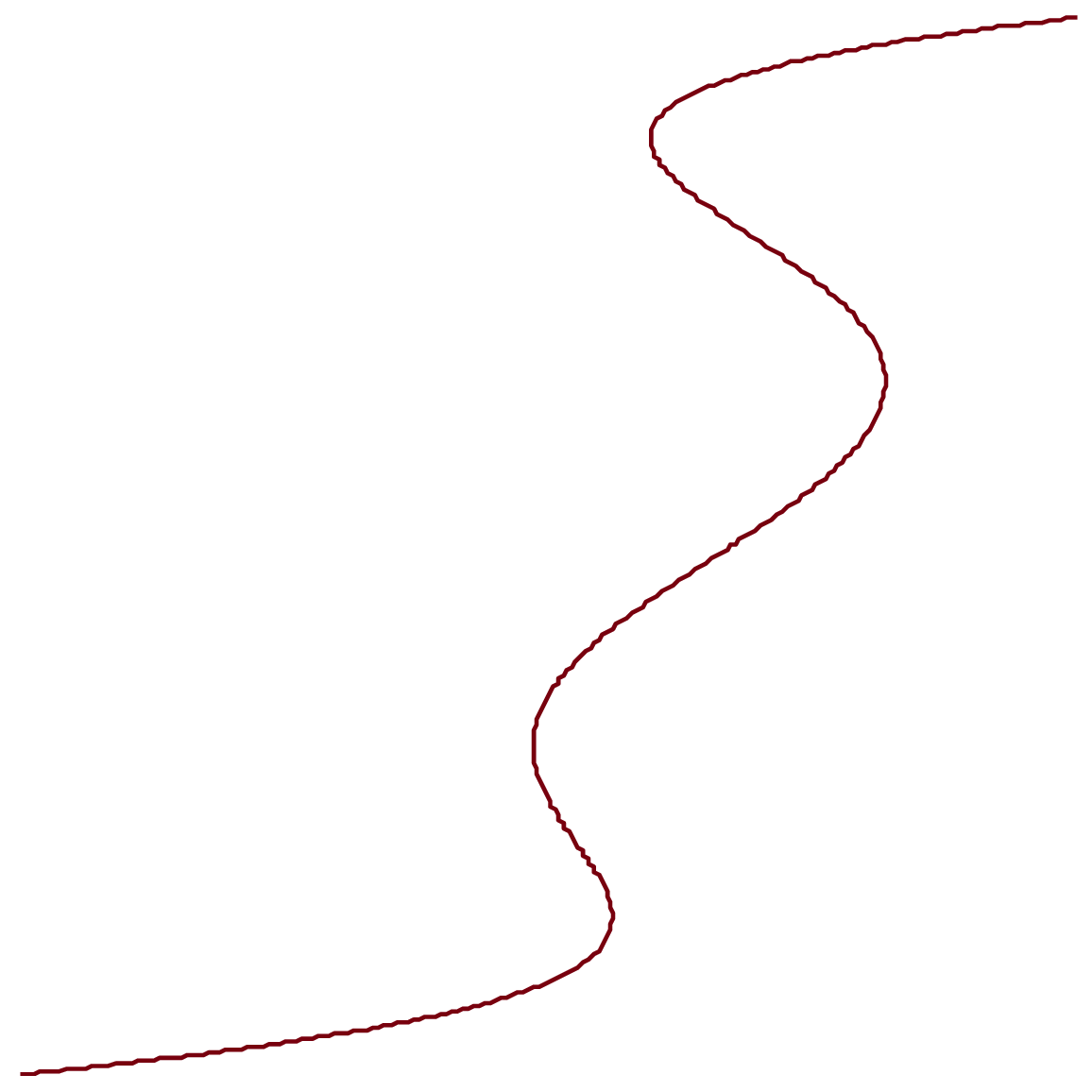}}
\subfigure[\label{e}]{\includegraphics[width=.26\columnwidth,height=.2\columnwidth]{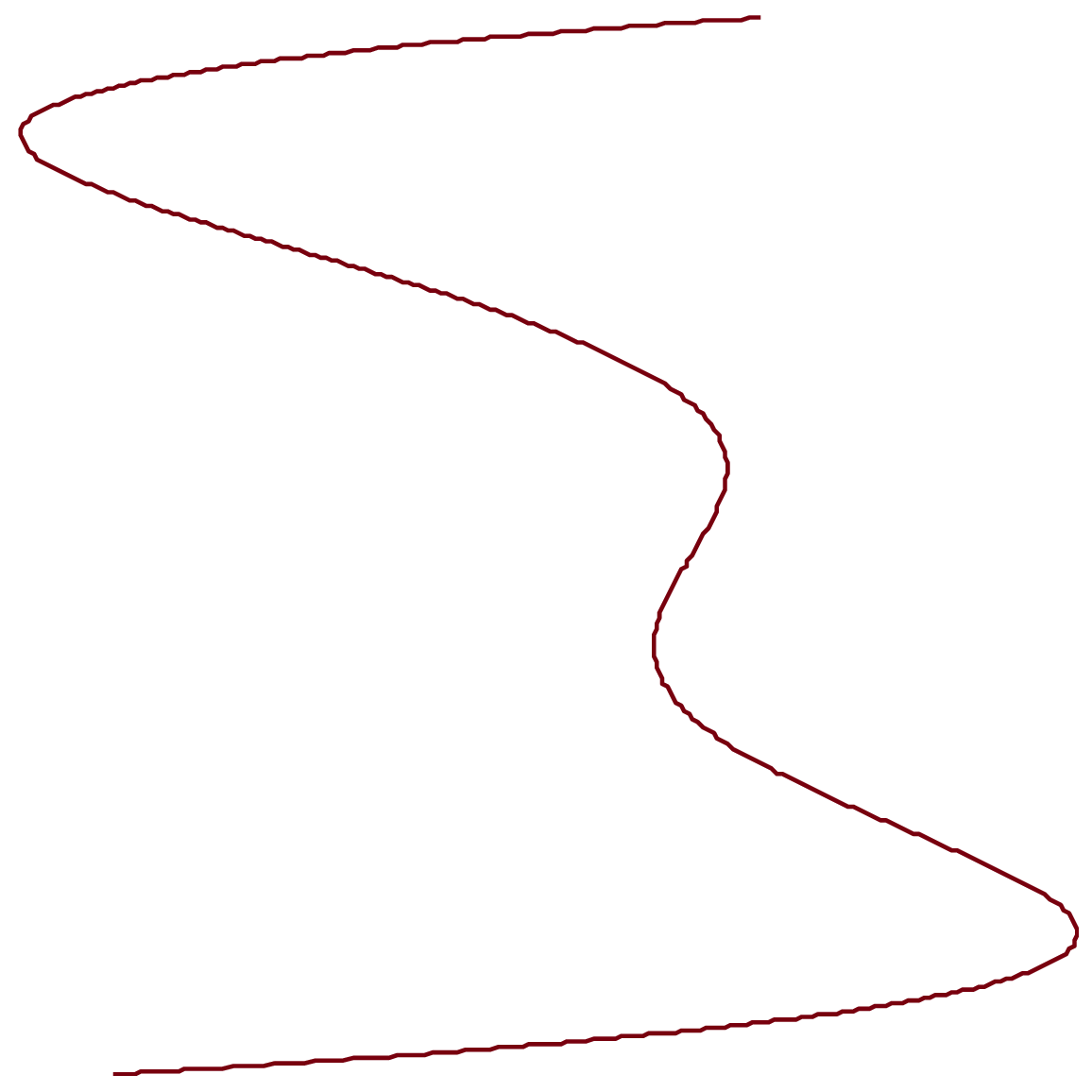}}
\subfigure[\label{f}]{\includegraphics[width=.26\columnwidth,height=.2\columnwidth]{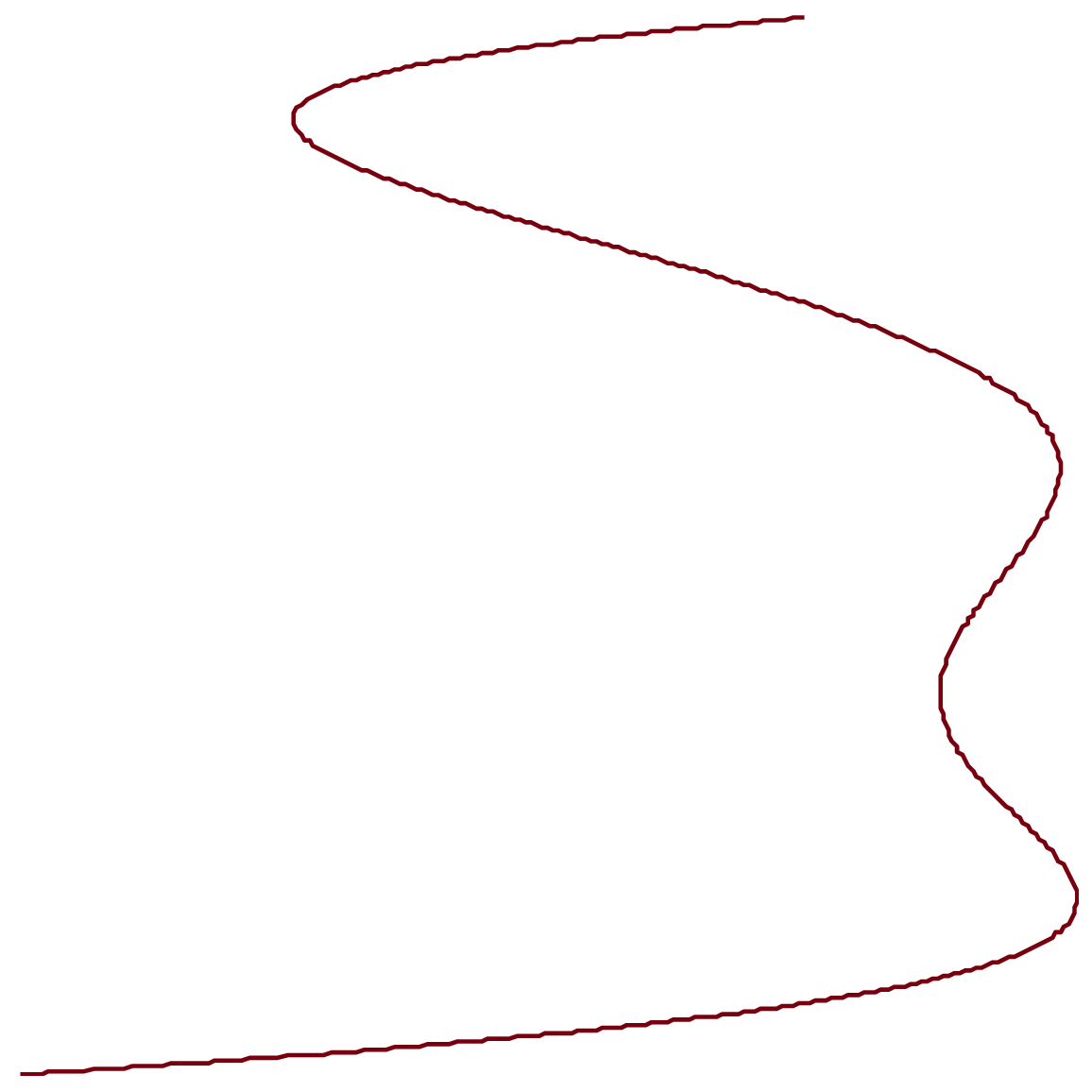}}
\subfigure[\label{g}]{\includegraphics[width=.26\columnwidth,height=.2\columnwidth]{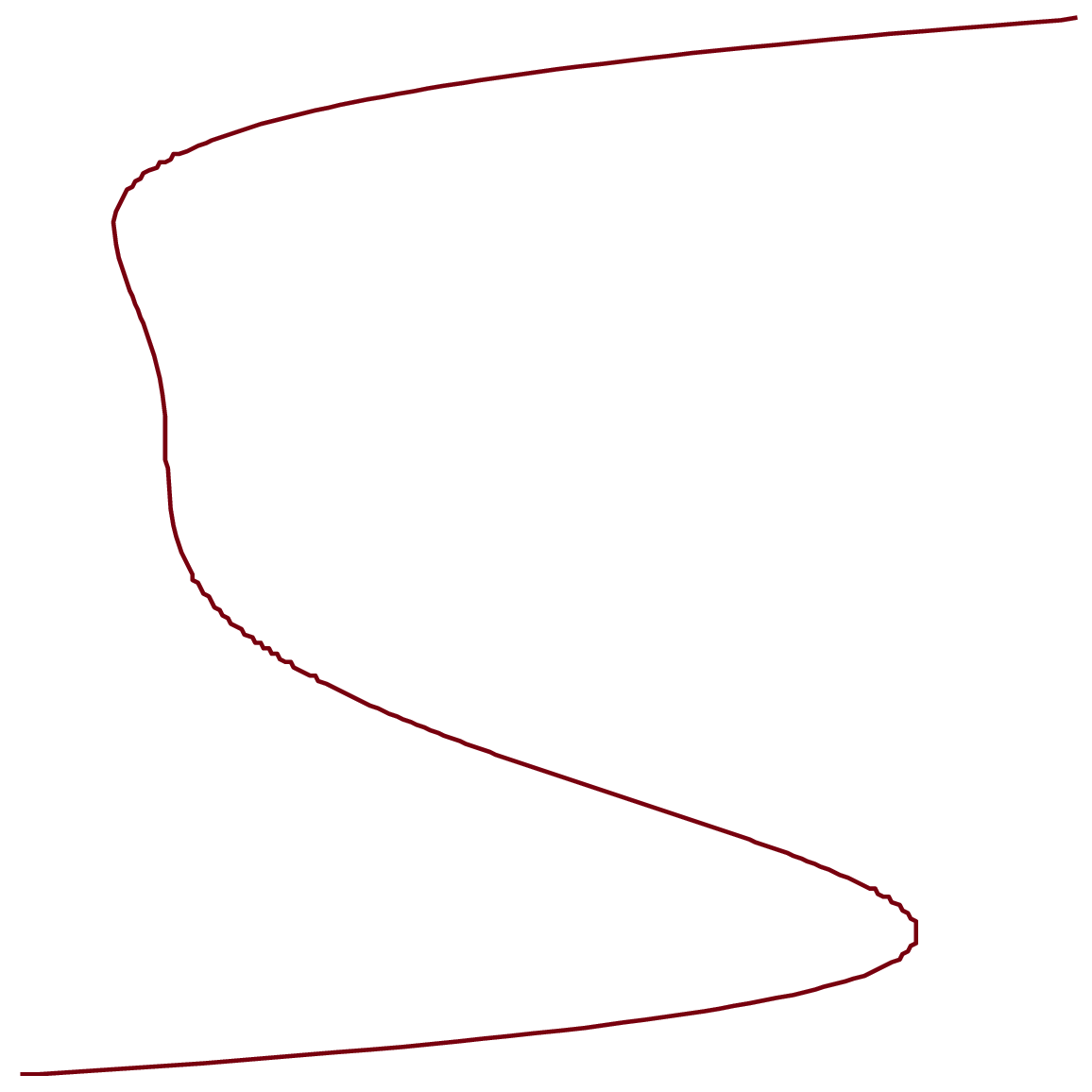}}
\subfigure[\label{h}]{\includegraphics[width=.26\columnwidth,height=.2\columnwidth]{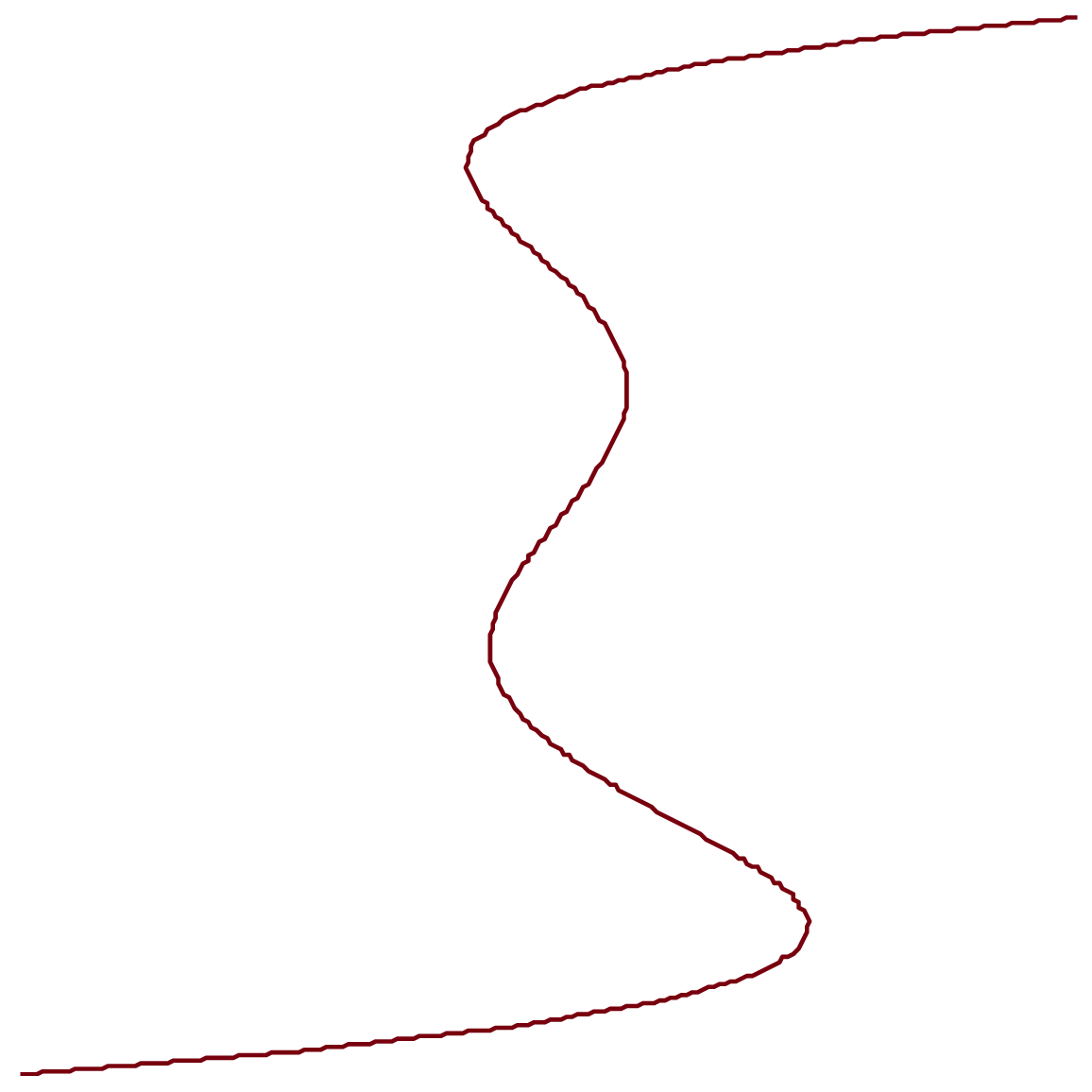}}
\subfigure[\label{i}]{\includegraphics[width=.26\columnwidth,height=.2\columnwidth]{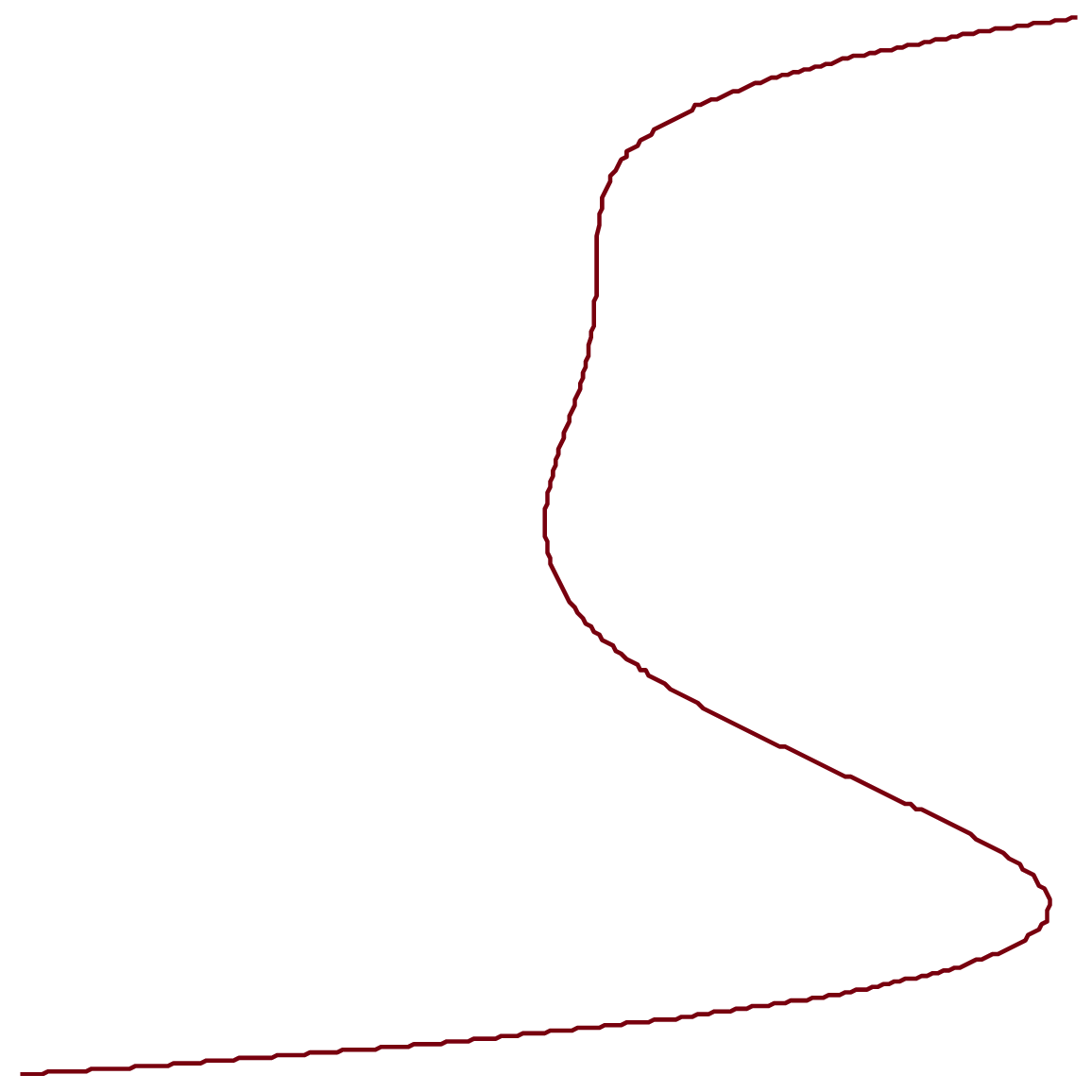}}
\caption{\small Persistent bifurcation diagrams associated with \(x^5-\lambda+\alpha_1 x+\alpha_2 x^2+\alpha_3 x^3\). }\label{Fig2}
\end{center}
\end{figure}
\end{example}

\subsection{Singular boundary conditions}\label{SecNonPersistent}

Extra sources of non-persistent is caused by singular boundary conditions of a parametric scalar map restricted to a bounded domain.
Let \(W\subset\mathbb{R}^{m}\) be a closed disk and \(U, L\subset\mathbb{R}\) be two closed intervals. Next, consider
\begin{equation*}
F(x,\lambda,\alpha)=0
\end{equation*}
where \(\alpha \in W\)  and \((x, \lambda)\in U \times L\); see \cite[Pages 154-158]{GolubitskySchaefferBook}. The new
non-persistent sources are defined by
\bas
  \mathscr{L}_{C}&:=& \{\alpha\in W\mid F(x,\lambda,\alpha)=0 \hbox{ for some }(x, \lambda)\in \partial U\times \partial L\},\\
  \mathscr{L}_{SH}&:=& \{\alpha\in W\mid F=F_{x}=0 \hbox{ for some }(x, \lambda)\in \partial U\times L\},\\
  \mathscr{L}_{SV}&:=& \{\alpha\in W\mid F=F_{x}=0 \hbox{ for some }(x, \lambda)\in U\times \partial L\},\\
  \mathscr{L}_{T}&:=& \{\alpha\in W\mid F=F_{\lambda}=0 \hbox{ for some }(x, \lambda)\in \partial U\times L\},
  \eas
  \bas
  \mathscr{G}_{1}&:=& \{\alpha\in W\mid F=0 \hbox{ at } (x_{0},\lambda,\alpha) \hbox{ for some }
  (x_{0},\lambda)\in \partial U \times L,\\
  && x_{0}\neq x \hbox{ and }
   F=F_{x}=0 \hbox{ at } (x, \lambda, \alpha) \hbox{ for some } (x, \lambda) \in U \times L\},\\
  \mathscr{G}_{2}&:=& \{\alpha\in W\mid\exists(x_{1}, \lambda), (x_{2}, \lambda)\in \partial U\times L, x_{1}\neq x_{2}\hbox{ s.t } F =0 \hbox{ at }\\
  && (x_{i}, \lambda, \alpha) \hbox{ for } i=1, 2\},
\eas
\bas
\mathscr{L_{B}}&:=& \{\alpha\in W\mid F=F_{x}=F_{\lambda}=0 \hbox{ at }  (x,\lambda,\alpha) \hbox{ for some }\\
  && (x, \lambda)\in U\times L\},\\
 \mathscr{L_{H}}&:=& \{\alpha\in W\mid F=F_{x}=F_{xx}=0 \hbox{ at } (x, \lambda, \alpha) \hbox{ for some }\\
 && (x, \lambda)\in U\times L\},\\
 \mathscr{G}_{D}&:=& \{\alpha\in W\mid\exists(x_{1}, \lambda), (x_{2}, \lambda)\in U\times L, x_{1}\neq x_{2}\hbox{ s.t } F= F_{x}=0 \hbox{ at }\\
  && (x_{i}, \lambda, \alpha) \hbox{ for } i=1,2\}.
  \eas
In this case, the transition set is given by \(\Sigma:=\mathscr{L}\cup\mathscr{G}\), here
\bas
\mathscr{L}&:=&\mathscr{L_{B}}\cup \mathscr{L_{H}}\cup\mathscr{L}_{C}\cup\mathscr{L}_{SH}\cup\mathscr{L}_{SV}\cup\mathscr{L}_{T},\\
\mathscr{G}&:=& \mathscr{G}_{D}\cup\mathscr{G}_{1}\cup \mathscr{G}_{2}.
\eas
For a finite codimension singular germ, \(\Sigma\) is a hypersurface of codimension one and each two choices from a connected component in the complement of \(\Sigma\) are contact-equivalent. Therefore, we can classify the persistent bifurcation diagrams by merely choosing one representative parameter from each components of the complement set of \(\Sigma\) and plotting the associated bifurcation diagrams. The command \verb"NonPersistent" is designed for this purpose.\\

\begin{tabularx}{\textwidth}{l|X}
  \textbf{Command/default option} & \textbf{Description} \\
\hline
\verb"NonPersistent"(\(F\), \(\alpha\), \verb"Vars", \(U\), \(L\)) & This function computes transition set for \(F\) where bifurcation
     diagrams are limited on \(U\times L\). Here, \(U\) and \(L\) are only taken as closed intervals. Further, it plots the transition set.
\\\hline
Box of figures & It  plots transition set in \([-1, 1]\times [-1, 1]\) by default.
\end{tabularx}
\subsubsection{Options}
\begin{itemize}
\item \(V, W\); this option enforces that the computed transition set is plotted in \(V\times W\) instead of the default square \([-1, 1]\times [-1, 1]\).
\item \verb"Vertical" (\verb"Horizontal" is also similar); this assumes that the boundary conditions is \(U\times \mathbb{R},\) \ie there is only singular boundary conditions on vertical boundary lines.
\end{itemize}

\begin{example}
\verb"NonPersistent"(\(x^4-\lambda x+\alpha_1 x+\alpha_2 \lambda+\alpha_3 x^2, [\alpha_1, \alpha_2, \alpha_3], [x, \lambda], [-2,2], [1,3]\)) gives rise to
\bas
 \mathscr{L}_{C}&=& \lbrace (\alpha_1, \alpha_2, \alpha_3)\,|\, \alpha_1+18+4\alpha_3+\alpha_2=0, \alpha_1+22+4\alpha_3+3\alpha_2=0, \\ && \alpha_1+14+4\alpha_3+\alpha_2=0, \alpha_1+10+4\alpha_3+3\alpha_2=0\rbrace,\\
 \mathscr{L}_{SH}&=& \lbrace (\alpha_1, \alpha_2, \alpha_3)\,|\, 4\alpha_2\alpha_3-\alpha_1+32\alpha_2+4\alpha_3+48=0, 4\alpha_2\alpha_3+\alpha_1\\&&+32\alpha_2-4\alpha_3-48=0\rbrace,\\
  \mathscr{L}_{SV}&=&\lbrace (\alpha_1, \alpha_2, \alpha_3)\,|\,  16\alpha_1\alpha_3^4+16\alpha_2\alpha_3^4-128\alpha_1^2\alpha_3^2-256\alpha_1\alpha_2\alpha_3^2
  -128\alpha_2^2\alpha_3^2\\&&+256\alpha_1^3+768\alpha_1^2\alpha_2+768\alpha_1\alpha_2^2
  +256\alpha_2^3-4\alpha_3^3+144\alpha_1\alpha_3+144\alpha_2\alpha_3-27=0,\\&& 16\alpha_1\alpha_3^4+48\alpha_2\alpha_3^4-128\alpha_1^2\alpha_3^2-768\alpha_1\alpha_2\alpha_3^2
  -1152\alpha_2^2\alpha_3^2+256\alpha_1^3+2304\alpha_1^2\alpha_2\\&&+6912\alpha_1\alpha_2^2
  +6912\alpha_2^3-36\alpha_3^3+1296\alpha_1\alpha_3+3888\alpha_2\alpha_3-2187=0 \rbrace,\\
  \mathscr{L}_{T}&=& \lbrace (\alpha_1, \alpha_2, \alpha_3)\,|\, \alpha_2+2=0, 4\alpha_3+16+\alpha_1=0, \alpha_2-2=0, 4\alpha_3+16+\alpha_1=0\rbrace,\\
   \mathscr{G}_{1}&=& \lbrace (\alpha_1, \alpha_2, \alpha_3)\,|\,4\alpha_2^2\alpha_3^3+80\alpha_2^2\alpha_3^2+16\alpha_2\alpha_3^3-72\alpha_1\alpha_2\alpha_3
   +512\alpha_2^2\alpha_3\\&&+32\alpha_2\alpha_3^2+16\alpha_3^3+27\alpha_1^2
   -320\alpha_1\alpha_2+72\alpha_1\alpha_3+1024\alpha_2^2-384\alpha_2\alpha_3
   \\&&+176\alpha_3^2+224\alpha_1-1024\alpha_2+640\alpha_3+768, 4\alpha_2^2\alpha_3^3+80\alpha_2^2\alpha_3^2-16\alpha_2\alpha_3^3\\&&+72\alpha_1\alpha_2\alpha_3
   +512\alpha_2^2\alpha_3-32\alpha_2\alpha_3^2+16\alpha_3^3+27\alpha_1^2+320\alpha_1\alpha_2
   +72\alpha_1\alpha_3\\&&+1024\alpha_2^2+384\alpha_2\alpha_3+176\alpha_3^2+224\alpha_1
   +1024\alpha_2+640\alpha_3+768=0
    \rbrace,\\
    \mathscr{G}_{2}&=& \lbrace (\alpha_1, \alpha_2, \alpha_3)\,|\,4\alpha_3+16+\alpha_1=0    \rbrace,\\
    \mathscr{L_{B}}&=& \lbrace (\alpha_1, \alpha_2, \alpha_3)\,|\,  \alpha_2^4+\alpha_2^2 \alpha_3+\alpha_1=0\rbrace,\\
     \mathscr{L_{H}}&=& \lbrace (\alpha_1, \alpha_2, \alpha_3)\,|\, 128\alpha_2^2\alpha_3^3+3\alpha_3^4+72\alpha_1\alpha_3^2+432\alpha_1^2=0  \rbrace,\\
     \mathscr{G}_{D}&=& \lbrace (\alpha_1, \alpha_2, \alpha_3)\,|\, \alpha_3^2-4\alpha_1=0 \rbrace,\\
\eas
and generates Figure \ref{Fig3}.
\begin{figure}[h]
\begin{center}
\subfigure[\label{1}]{\includegraphics[width=.40\columnwidth,height=.4\columnwidth]{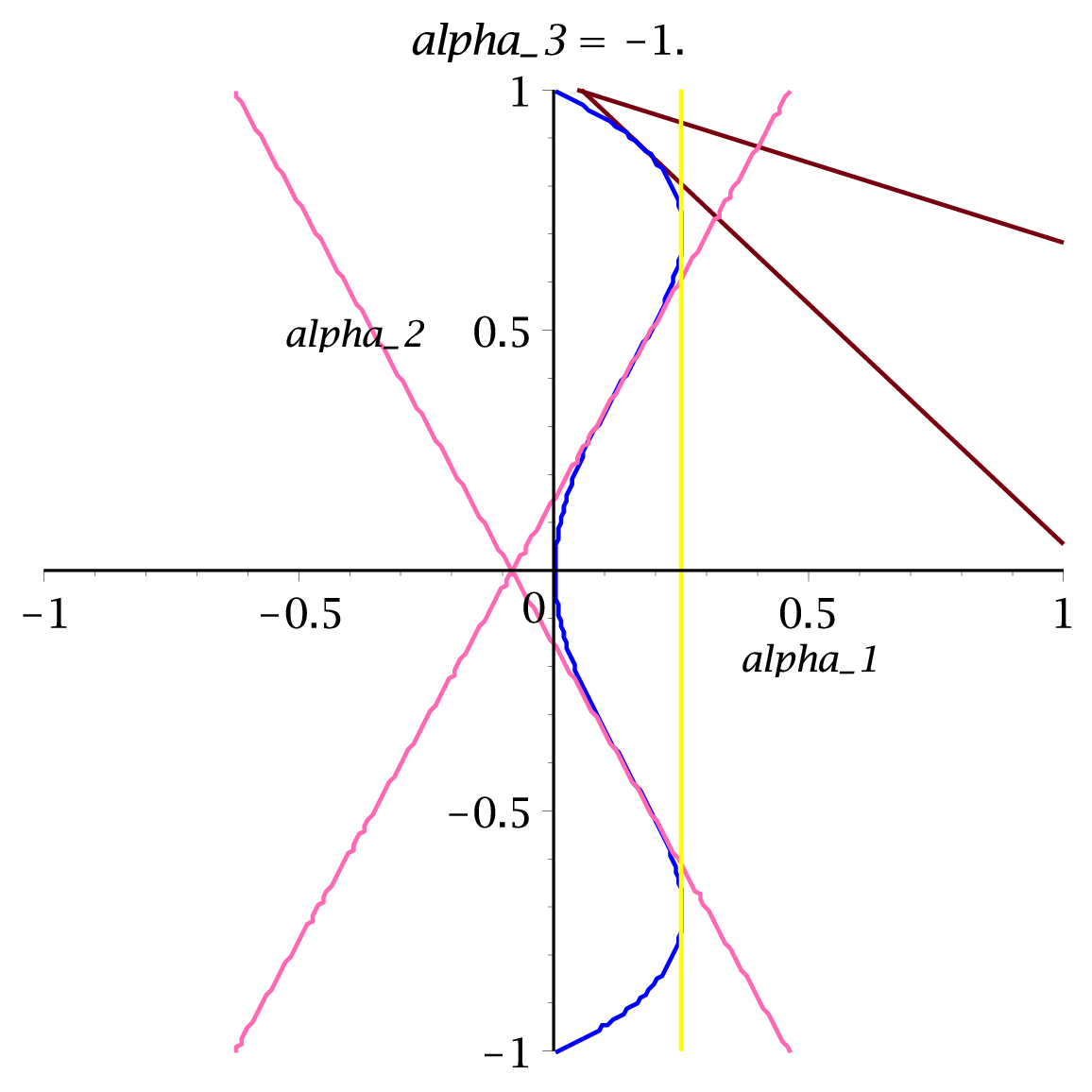}}
\caption{\small Transition set}\label{Fig3}
\end{center}
\end{figure}
\end{example}

\part{Tools from algebraic geometry }

In this section we describe how to compute some tools from computational algebraic geometry; see \cite{GazorKazemi} for more information.

\section{Multiplication Matrix}\label{SecMultMatrix}

Let \(\mathcal{R}\) be either of the rings of germs \(\E,\) \(\mathbb{R}[[x, \lambda]],\) \(\mathbb{R}[x, \lambda]_{\langle x, \lambda\rangle}\) or \(R[x, \lambda]\); see \cite{GazorKazemi}. Now we describe how to compute the multiplication Matrix defined by
\be\label{Mult}
\varphi_{u, J}: \frac{\mathcal{R}}{J}\rightarrow \frac{\mathcal{R}}{J}, \quad \varphi_{u, J}(f+J):= uf+ J,
\ee where \(J\) is an ideal generated by a finite set \(A\subset \mathcal{R},\) \ie \(J:=\langle A\rangle_\mathcal{R}\), and \(u\) is a monomial; also see \cite[Equation 3.4]{GazorKazemi}.\\

\begin{tabularx}{\textwidth}{l|X}
  \textbf{Command/option} & \textbf{Description} \\
      \hline
  \verb"MultMatrix"(\(A\), \(u\), \verb"Vars") & This function derives \(\varphi_{u, J}\) where \(u\) is a monomial, \(\varphi_{u, J}\) is defined by Equation \eqref{Mult}. \\\hline
Default computational ring & the fractional germ ring. \\
      \hline
Truncation degree & When the input set of germs \(A\) only includes polynomials, \verb"MultMatrix"(\(A\), \(u\), \verb"Vars") does not need truncation degree. However, for non-polynomial input germs, a truncation degree \(k\) needs to be included. \\
\end{tabularx}

\subsection{Options}
\begin{itemize}
  \item \(k\); determines the truncation degree. The user is advised to use the command \verb"Verify" to find an appropriate truncation degree \(k.\)
  \item Computational ring: \verb"Fractional", \verb"Formal", \verb"SmoothGerms"; \verb"Polynomial"; the command uses either the rings of fractional germs, formal power series or ring of smooth germs. The command \verb"Verify" is an appropriate tool to find/verify the appropriate computational ring.
\end{itemize}

\begin{example}
The command \verb"MultMatrix"(\([x^6+\frac{12}{27}x^10\lambda^9, \frac{5}{3}x^{5}+\lambda
\sin(x^3), \lambda^2-\frac{2}{3}(1-\exp(x^5))], x, [x, \lambda] , 6, \verb"Formal"\)) gives rise to

\bes
\left( {\begin{array}{ccccccccc}
0 & 0 & 0 & 0 & 0 & 0 & 0 & 0 & 0\\
0 & 0 & 0 & 0 & 0 & 0 & 0 & 0 & 0\\
1 & 0 & 0 & 0 & 0 & 0 & 0 & 0 & 0\\
0 & 0 & 1 & 0 & 0 & 0 & 0 & 0 & 0\\
0 & 0 & 0 & 1 & 0 & 0 & 0 & 0 & 0\\
0 & 0 & 0 & 0 & 1 & 0 & 0 & 0 & 0\\
0 & 0 & 0 & 0 & 0 & 1 & 0 & 0 & -\frac{5}{3}\\
0 & 1 & 0 & 0 & 0 & 0 & 0 & 0 & 0\\
0 & 0 & 0 & 0 & 0 & 0 & 0 & 1 & 0\\
\end{array}}
\right).
\ees

\end{example}

\section{Divisions of germs }\label{SecDivision}

The following table describes how to use the command \verb"Division" to divide a germ \(g\) by a set of germs \(G:= \{f_i| i=1, \ldots n\}\) where all these germs are in terms of the variables in \verb"Vars". Note that the ordering in the list of variables \verb"Vars" is important and determines how anti-lexicographical ordering is defined.\\

\begin{tabularx}{\textwidth}{l|X}
  \textbf{Command} & \textbf{Description} \\
\hline
\verb"Division"(\(g\), \(G\), \verb"Vars") & This divides the germ \(g\) by the germs in \(G\) using anti-lexicographical ordering.  \\
\hline
\end{tabularx}

\subsection{Options}
\begin{itemize}
\item \verb"Formal"; \verb"SmoothGerms"; \verb"Polynomial"; \verb"Fractional"; this determines the computational germ ring. In order to verify/check the computational ring for the division, one is advised to use the command \verb"Verify"(G, \verb"Ideal", \verb"Vars") and find the permissible computational ring.
\item \(k\); this option enforces the computations modulo degrees higher than or equal to \(k+1.\) The command \verb"Verify"(G, \verb"Ideal", \verb"Vars") also suggests an optimal permissible truncation degree \(k\).
\end{itemize}
\begin{example}
For an example we use \verb"Division"(\(\sin(x^7)-1, [x^5+x^6\exp(\lambda), xy^3-\frac{2}{7}x\lambda^6-x^7, \lambda\cos(x^7)],\) \([x, \lambda], 8,\verb"SmoothGerms"\)). This returns \bes -1\ees as the remainder of \bes \sin(x^7)-1\ees divided by
\bes \{x^5+x^6\exp(\lambda), xy^3-\frac{2}{7}x\lambda^6-x^7, \lambda\cos(x^7)\}\ees using the anti-lexicographical ordering with \(x\succ\lambda\).
\end{example}

\section{Standard basis} \label{SecStandardBasis}

Now we describe how to compute a standard basis for a set of germs in either of the following local rings: fractional germs \(\mathbb{R}[x, \lambda]_{\langle x, \lambda\rangle}\), formal power series \(\mathbb{R}[[x, \lambda]],\) and ring of smooth germs \(\E\); see \cite{GazorKazemi}. The command \verb"StandardBasis" follows the table and options described below. Note that \verb"Vars" denotes an order list of variables and the germs in \(G\) are in terms of the variables in \verb"Vars". The anti-lexicographical ordering is here used while it is determined by the ordering of variables in \verb"Vars".\\

\begin{tabularx}{\textwidth}{l|X}
  \textbf{Command} & \textbf{Description} \\
\hline
\verb"StandardBasis"(\(G\), \verb"Vars") & computes the standard basis of the polynomial germs in \(G\) in the fractional germ ring. \\
\hline
Default computational ring & fractional germ ring. \\
\hline
Default input & the default germs in \(G\) must be polynomial germs. For the cases of non-polynomials, it needs a truncation degree \(k.\)
\end{tabularx}

\subsection{Options}

\begin{itemize}
\item \(k\); this determines the truncation degree.
\item Consider the cases that the option \(k\) is used. A warning note is given when the truncation degree \(k\) is not sufficiently high to guarantee that
the output is correct. A warning note is given like ``The truncation degree is not sufficiently high and thus, the following results might be wrong.'' In this case, for an appropriate truncation degree \(l\) is given.
\item \verb"Formal"; \verb"SmoothGerms"; \verb"Fractional"; this determines the computational germ ring and computes the standard basis accordingly.
\end{itemize}
\begin{example}
For an example, the command
\verb"StandardBasis"(\([x^5+x^2\sin(\lambda+x)+\lambda^2, x^3\lambda^2+\cos(\lambda)x, \lambda^6+x^4-\lambda x], [x, \lambda], 7, \verb"SmoothGerms"\))
computes the standard basis of the set of germs
\bes \{x^5+x^2\sin(\lambda+x)+\lambda^2, x^3\lambda^2+\cos(\lambda)x, \lambda^6+x^4-\lambda x\}\ees as
\(\{x, \lambda^2\}.\)
\end{example}

\section{Colon ideals}\label{SecColonIdeal}

The colon ideal \(I: \langle g\rangle \) refers to the ideal defined by
\bes
I: \langle g\rangle_\E= \{f\in \E: f \langle g\rangle_\E\subseteq I\}.
\ees Using the arguments on \cite[Page 22]{GazorKazemi}, we have \(I: \langle g\rangle_\E= \langle \frac{h_i}{g}| i=1, \ldots, n\rangle_\E,\) where \(h_1, h_2, \ldots, h_n\) is a standard basis for the ideal \(I\cap \langle g\rangle_{\mathbb{R}[x, \lambda]}.\) The command \verb"ColonIdeal" follows this. \\

\begin{tabularx}{\textwidth}{l|X}
  \textbf{Command} & \textbf{Description} \\
\hline
\verb"ColonIdeal"(\(I\), \(g\))& computes the colon ideal \(I:g\) for \(g\) in \(\E\).
\end{tabularx}
\begin{example}
\verb"ColonIdeal"(\([x^7+\lambda x^3-\lambda^2 x, \lambda x^6+\lambda^2 x^2-\lambda^3, x^3\lambda+x],\lambda\)) leads to the ideal generated by
\bes
[x(\lambda x^2+1), \lambda^2 x^2-x^4-\lambda^3, \lambda^4+\lambda^2-x^2, x(x^4+\lambda^3+\lambda)].
\ees
\end{example}

\section{Complement spaces }\label{SecNormalset}

For computing universal unfolding of a singular germ \(g\), we need to compute a basis for a complement vector space for the tangent space \(T(g)\) associated with \(g.\) This is equivalent to computing a basis for the quotient space \(\E/T(g).\) More generally, the command
\verb"Normalset"(\(I\)) computes a monomial basis for \(\E/I\), when \(I\) is either an ideal or a vector space with finite codimension in the local ring \(\E\).\\

\begin{tabularx}{\textwidth}{l|X}
  \textbf{Command} & \textbf{Description} \\
\hline
\verb"Normalset"(\(A\)) & computes a monomial basis for \(\E/I\), when \(A\) is a list of germs generating an ideal \(I:=\langle A\rangle\).
\end{tabularx}

\begin{example}
For example the command \verb"Normalset"(\([x^6+\lambda x^4+\lambda^2 x, \lambda x^5+\lambda^2 x^3+\lambda^3, 5x^6+3\lambda x^4, 5\lambda x^4+3\lambda^2 x^2, -3x^7-3\lambda x^5-\frac{25}{3}x^6]\)) returns
\bes
[1, \lambda, x, \lambda^2, x^2, x^3, x^4, x^5, \lambda x, \lambda x^3, x^2\lambda]
\ees as a list of monomials for the complement of \(I\) in \(\E.\)
\end{example}

\part{Technical objects in singularity theory}

A user who is not interested in technical details and their commands may simply skip this section.

\section{Intrinsic ideals }\label{SecIntrinsic}

In this section we describe how to compute the intrinsic part of an ideal or a vector space. Let \(A\) and \(B\) be two lists indicating the generators of an ideal and a vector space, respectively. We intend to compute the maximal intrinsic ideal contained in \(\langle A\rangle_\E+\langle A\rangle_{\mathbb{R}}.\)\\

\begin{tabularx}{\textwidth}{l|X}
  \textbf{Command/option} & \textbf{Description} \\
\hline
\verb"Intrinsic"(\(A\), \verb"Vars")  & computes intrinsic part of the ideal generated by \(A\). It remarks when the ideal is of infinite codimension. \\
 \hline
\verb" Intrinsic"(\(A, B\), \verb"Vars") & computes intrinsic part of a vector space given by \(\langle A\rangle_\E + \langle B\rangle_{\mathbb{R}}.\) It remarks when the vector space is of infinite codimension. \\
 \hline
 Computational ring &The default computational ring is the ring of fractional germs. \\
 \hline
Verify/Warning/Suggestion& By default, \verb" Intrinsic" checks and verifies whether the fractional germ ring is sufficient for computation of the intrinsic part of the ideal or vector space spanned by \(A\) (and \(B\)). If fractional germ ring is not sufficient, it gives a warning note of possible errors and a suggestion to circumvent the problem. \verb" Intrinsic" remarks when the problem might be of infinite codimension. Despite possible warning errors, \verb" Intrinsic" computes the intrinsic part using the fractional germ ring.
\end{tabularx}

\subsection{Options}
\begin{itemize}
\item \(k\); this option enforces that the computations are performed modulo degree \(k.\)
\item \verb"Formal"; \verb"SmoothGerms"; \verb"Polynomial"; \verb"Fractional"; this determines the computational germ ring. It checks and verifies the computations. It gives warning notes of possible errors and alternative suggestions when it finds them necessary.
\end{itemize}

\begin{example} For an example
\verb"Intrinsic"\(([x^3\lambda+\lambda^2, 3x^3\lambda, 3x^2\lambda^2],\) \([x, \lambda]\), \(\verb"InfCodim"\), \(\verb"Fractional")\)  leads to
\bes
M^3\langle\lambda\rangle+\langle\lambda^2\rangle.
\ees
As for a second example \verb"Intrinsic"(\([x^5+\lambda x^3+\lambda^2, 5x^5+3x^3\lambda, 5x^4\lambda+3x^2\lambda^2], [\lambda x^3+2\lambda^2, x^3+2\lambda, x^4+\frac{3}{5}\lambda x^2, \lambda^2, x^5], [x, \lambda]\)) results in
\bes
M^5+M^3\langle \lambda\rangle+\langle\lambda^2\rangle.
\ees

\end{example}

\section{Algebraic objects }\label{SecAlgObjects}

Singularity theory defines and uses many algebraic objects in the bifurcation analysis of zeros of smooth germs. These include
restricted tangent space \(RT\), tangent space \(T\), high order term ideal \(\p\), smallest intrinsic ideal associated with a singular germ \(\s\), a basis for complement of the tangent space \(\E/T\), and low order terms \(\s^{\perp}\); see \cite{GazorKazemi,GolubitskySchaefferBook}. These can be computed in {\tt Singularity} using  the command \verb"AlgObjects" as well as the individual commands \verb"RT", \verb"T", \verb"P", and \verb"TangentPerp". The individual commands \verb"RT", \verb"T", and \verb"P" have the same default and non-default options as \verb"AlgObjects" has as follows.\\

\begin{tabularx}{\textwidth}{l|X}
  \textbf{Command/option} & \textbf{Description} \\
\hline
 \verb"AlgObjects"(\(g\), \verb"Vars")  & This function computes \(RT\), \(T\), \(\p\), \(\E/T\), \(\s\), \(\s^{\perp}\) and intrinsic generators of \(\s\) for given \(g\).
\\
\hline
\verb"RT"(g, \verb"Vars") & This derives the restricted tangent space associated with a scalar smooth germ \(g\).\\
 \hline
\verb"T"(g, \verb"Vars") & This command provides a nice representation of the tangent space associated with the singular smooth germ \(g.\) The representation uses intrinsic ideal representation as for the intrinsic part of \(T(g)\).
\\
 \hline
\verb"p"(g, \verb"Vars") & This computes the high order term ideal associated with the germ \(g\). \\
\hline
\verb"TangentPerp"(g, \verb"Vars")& This first computes \(T(g)\), \ie the tangent space of the germ \(g,\) and then returns a monomial basis for the complement space of \(T(g)\).\\
\hline
\verb"S"(g, \verb"Vars")& Computes the smallest intrinsic ideal containing the germ \(g\).\\
\hline
\verb"SPerp"(g, \verb"Vars")& This derives a set of monomials of low order terms for the germ \(g\).\\
\hline
\verb"IntrinsicGen"(g, \verb"Vars")& This derives the intrinsic generators of \(\s(g)\) that determine the nonzero conditions for recognition problem for normal forms.
\\
 \hline
Computational ring & The default computational ring is the ring of fractional germs. \\
 \hline
Default degree \(k\) & For non-polynomial input germs, it computes the least degree \(k\) so that truncations at degree \(k\) is permissible. Next, the germ \(g\) is truncated and all algebraic objects are computed modulo degrees higher than or equal to \(k+1\).
\end{tabularx}

\subsection{Options}
\begin{itemize}
\item \(k\); this option enforces the computations modulo degree \(k+1.\)
\item \verb"Formal"; \verb"SmoothGerms"; \verb"Polynomial"; \verb"Fractional"; this determines the computational germ ring. It checks and verifies/gives warning notes of possible errors.
\end{itemize}

\begin{example} Now we present three examples of singular germs of high codimension; see and compare these examples with the examples in \cite[Page 4]{Keyfitz}. For example we consider a codimension 10 singularity and use \verb"AlgObjects"(\(x^5+x^3\lambda^2+\lambda^3, [x, \lambda],  \verb"Fractional"\)). It gives
\bas
\p &=& M^{6}+M\langle \lambda^3 \rangle,\\
RT &=& M^{6}+M\langle \lambda^3 \rangle+\langle x^4\lambda, 3\lambda^2x^3+5x^5, \lambda^2x^3+x^5+\lambda^3 \rangle,\\
T&=& M^{5}+\langle \lambda^3\rangle+\lbrace \frac{3}{5}\lambda^2x^2+x^4, x^3\lambda+\frac{3}{2}\lambda^2\rbrace,\\
\E/T &=&\lbrace 1, \lambda, x, \lambda^2, x^2, x^3, \lambda x^2, \lambda^2 x, \lambda^2x^2, x\lambda\rbrace,\\
\s &=& M^{5}+\langle \lambda^3\rangle,\\
\s^{\perp} &=&\lbrace 1, \lambda, x, \lambda^2, x^2, x^3, x^4, \lambda x^2, \lambda^2 x, \lambda^2 x^2, x\lambda, x^3\lambda\rbrace,\\
\hbox{ intrinsic generators } &=& \lbrace x^5, \lambda^3\rbrace.
\eas
A codimension 20 singularity:
\verb"TangentPerp"(\(x^8+\sin(\lambda^3), [x, \lambda], 9\)) derives the following
\bas
T&=& M^{9}+\langle \lambda^3\rangle+\lbrace x^7, x^8, \lambda x^7, -2\lambda^2\rbrace,\\
\hbox{ TangentPerp }&=&\lbrace 1, \lambda, x, x^2, x^3, x^4, x^5, x^6, \lambda x, \lambda x^2, \lambda x^3, \lambda x^6,\\&& \lambda^2 x, \lambda^2 x^2, \lambda^2 x^5, x^3\lambda^2, x^4\lambda, x^4\lambda^2, x^5\lambda, x^6\lambda^2\rbrace.
\eas

Now we use the command \verb"IntrinsicGen" for an example of codimension 13. \verb"IntrinsicGen"(\(\lambda x^8+x^7-\lambda^3 x^2-\lambda^2 x, [x, \lambda]\)) leads to
\bas
\hbox{ intrinsic generators } &=& \lbrace x^7, x\lambda^2 \rbrace.
\eas
\end{example}
\section*{\bf Acknowledgments}

We thank Professor Erik Postma at Maplesoft for his time and fruitful comments for
improvement of our library.

\end{document}